\documentclass[11pt,notitlepage,twoside,a4paper]{amsart}
 \usepackage{amsfonts}

\usepackage{amsmath,amssymb,enumerate}

\usepackage{epsfig,fancyhdr,color}%,showkeys,amsmidx

\usepackage{amssymb}
\usepackage{amsmath,amsthm}
\usepackage{latexsym}
\usepackage{amscd}
\usepackage{psfrag}
\usepackage{graphicx}
\usepackage[latin1]{inputenc}
\usepackage[all]{xy}
\usepackage[mathcal]{eucal}

\definecolor{NoteColor}{rgb}{1,0,0}

\renewcommand{\textsc}{\textcolor{red}}

\newtheorem*{theorem 1}{\rm\bf Proposition 1}
\newtheorem*{theorem 2}{\rm\bf Proposition 2}

\theoremstyle{definition}

\theoremstyle{remark}

\def\interieur#1{\mathord{\mathop{\kern 0pt #1}\limits^\circ}}

\title[Quasiconformal mappings]{Quasiconformal mappings,  
\\ 
 from Ptolemy's geography to the work of Teichm\"uller}

\author[Athanase Papadopoulos]{Athanase Papadopoulos}

\thanks{The work was completed during a stay of the author at the the Max-Plank-Institut f\"ur Mathematik (Bonn). The author thanks this institute for its hospitality.}

 \address{Athanase Papadopoulos, Max-Plank-Institut f\"ur Mathematik, Vivatsgasse 7, 53111 Bonn, Germany, and Institut de Recherche Math{\'e}matique Avanc\'ee,
Universit{\'e} de Strasbourg and CNRS,
7 rue Ren\'e Descartes,
 67084 Strasbourg Cedex, France} \email{papadopoulos@math.u-strasbg.fr}

\date{\today}

\begin{document}

  \maketitle

  \hfill{\emph{Les hommes passent, mais les \oe uvres restent}.}
  
    \hfill{(Augustin Cauchy, \cite{Cauchy-Vie} p. 274)}

  \begin{abstract}
  
 The origin of quasiconformal mappings, like that of conformal mappings, can be traced back to old cartography where the basic problem was the search for mappings from the sphere onto the plane with minimal deviation from conformality,  subject to certain conditions which were made precise. In this paper, we survey the development of cartography, highlighting the main ideas that are related to quasiconformality. Some of these ideas were completely ignored in the previous historical surveys on quasiconformal mappings. We then survey early quasiconformal theory in the works of Gr\"otzsch, Lavrentieff, Ahlfors and Teichm\"uller, which are the 20th-century founders of the theory.

  The period we consider starts with Claudius Ptolemy (c. 100--170 A.D.) and ends with Oswald Teichm\"uller (1913--1943). We mention the works of several mathematicians-geographers done in this period, including Euler, Lagrange, Lambert, Gauss, Chebyshev, Darboux and others.  We highlight in particular  the work of Nicolas-Auguste Tissot  (1824--1897), a French mathematician and geographer who (according to our knowledge)  was the first to introduce the notion of a mapping which transforms infinitesimal circles into infinitesimal ellipses, studying parameters such as the ratio of the major to the minor axes of such an infinitesimal ellipse, its area divided by the area of the infinitesimal circle of which it is the image, and the inclination of its axis with respect to a fixed axis in the plane.
 We also give some information about the lives and works of Gr\"otzsch, Lavrentieff, Ahlfors and Teichm\"uller. The latter brought the theory of quasiconformal mappings to a high level of development. He used it in an essential way in his investigations of Riemann surfaces and their moduli and in function theory (in particular, in his work on the Bieberbach conjecture and the type problem). We survey in detail several of his results. We also discuss some aspects of his life and writings, explaining why his papers were not read and why some of his ideas are still unknown even to Teichm\"uller theorists. 
 
 The final version of this paper will appear in the book ``Uniformization, Riemann-Hilbert Correspondence, Calabi-Yau Manifolds, and Picard-Fuchs Equations" (ed. L. Ji and S.-T. Yau),  International Press and  Higher Education Press (2017). 

  \end{abstract}

\smaller{
\noindent AMS Mathematics Subject Classification:  01A55, 30C20, 53A05, 53A30, 91D20. 

\noindent Keywords:  Quasiconformal mapping, geographical map, sphere projection, Tissot indicatrix.
\larger}

 \tableofcontents
 
   \section{Introduction}

  \medskip

  \medskip

The roots of quasiconformal theory lie in geography, more precisely in the study of mappings from (subsets of)  the sphere to the Euclidean plane, and the attempts to find the ``best" such mappings. The word ``best" refers here to  mappings with the least ``distortion," a notion to be defined appropriately, involving distances, angles and areas. 
The period we cover starts with Greek antiquity and ends with the work of Teichm\"uller. The latter, during his short lifetime, brought the theory of quasiconformal mappings to matureness. He made an essential use of it in his work on Riemann's moduli problem \cite{T20}   \cite{T29} \cite{T24}, on the Bieberbach coefficient problem \cite{T14}, on value distribution theory \cite{T200} \cite{T-Einfache} and on the type problem \cite{T9}. The so-called Teichm\"uller Existence and Uniqueness Theorem concerns the existence and uniqueness of the best quasiconformal mapping between two arbitrary Riemann surfaces in a given homotopy class of homeomorphisms. This result is a wide generalization of a result of Gr\"otzsch on mappings between rectangles \cite{Gr1928} \cite{Gr}. Teichm\"uller proved this theorem in the two papers \cite{T20} and \cite{T29}. Motivated by his work on extremal quasiconformal mappings, he wrote a paper \cite{T23} on extremal problems from a very general broad point of view. In that paper, he made analogies between algebra, Galois theory, and function theory, based on the idea of extremal mappings.  He also wrote a paper on the solution of a specific problem on  quasiconformal mappings \cite{T31}, which has several natural generalizations.

Regarding Teichm\"uller's  work, Ahlfors, in his 1953 paper on the \emph{Development of the theory of conformal mapping and Riemann surfaces through a century}, written at the occasion of the hundredth anniversary of Riemann's inaugural dissertation \cite{Ahlfors-dev},  says:
\begin{quote}\small
In the premature death of Teichm\"uller, geometric function theory, like other branches of mathematics, suffered a grievous loss. He spotted the image of Gr\"otzsch's technique, and made numerous applications of it, which it would take me too long to list. Even more important, he made systematic use of extremal quasiconformal mappings, a concept that Gr\"otzsch had introduced in a very simple special case. Quasiconformal mappings are not only a valuable tool in questions connected with the type problem, but through the fundamental although difficult work of Teichm\"uller it has become clear that they are instrumental in the study of modules of closed surfaces. Incidentally, this study led to a better understanding of the role of quadratic differentials, which in somewhat mysterious fashion seem to enter in all extremal problems connected with conformal mapping.
\end{quote}

Talking also about Teichm\"uller, Lehto writes in \cite{Lehto1984} (1984), p. 210: 
                                           \begin{quote}\small
                                           In the late thirties, quasiconformal mappings rose to the forefront of complex analysis thanks to Teichm\"uller. Introducing novel ideas, Teichm\"uller showed the intimate interaction between quasiconformal mappings and Riemann surfaces.
                                           \end{quote}
                                           
                     In the present paper, while trying to convey the main ideas and mathematical concepts related to quasiconformal mappings, we make several digressions which will help the reader including these ideas in their proper context. The largest digression concerns the history of geographical maps. We tried to show the broadness of the subject and its importance in the origin of differential geometry.

\medskip

The plan of this paper is as follows:

\medskip

After the present introduction, the paper is divided into 8 sections. 

Section \ref{s:qc} is an exposition of the history of quasiconformal mappings, before the name was given, and their origin in geography.   We start in Greek antiquity, namely, in the work of Ptolemy on the drawing of geographical maps. Ptolemy quotes in his \emph{Geography} the works of several of his predecessors on the subject, and we shall mention some of them. We then survey  works of some major later mathematicians who contributed in various ways to the development of cartography. The names include Euler, Lagrange, Chebyshev, Gauss, Darboux and Beltrami. The long excursion that we take in geography and the drawing of geographical maps will also show that we find there the origin of several important ideas in differential geometry. Furthermore, the present article will show the close relationship between geometry and problems related to the measurement of the Earth, a relationship which was very strong in antiquity and which never lost its force, and which confirms the original meaning of the word \emph{geo-metria}\footnote{Proclus (5th c. A.D.), in his \emph{Commentary on Book I of Euclid's Elements}, considers that geometry was discovered by the Egyptians, for the needs of land surveying. Proclus' treatise is sometimes considered as the oldest essay on the history of geometry. It contains an overview of the subject up to the time of Euclid.} (earth-measurement).

\S \ref{s:Tissot} concerns the work of  Nicolas-Auguste Tissot. This work makes the link between old cartography and the modern notion of quasiconformal mappings. Tissot associated to a map $f:S_1\to S_2$ between surfaces the field of infinitesimal ellipses on $S_2$ that are the images by the differential of the map $f$ of a field of infinitesimal circles on $S_1$. He used the geometry of this field of ellipses as a measure of the distortion of the map.  At the same time, he proved several properties of the ``best map"  with respect to this distortion, making use of a pair of perpendicular foliations on each of the two surfaces that are preserved by the map. He showed the existence of these pairs of foliations, and their uniqueness in the case where the map is not conformal. All these ideas were reintroduced later on as essential tools in modern quasiconformal theory. Tissot's pair of foliations are an early version of the pair of measured foliations underlying the quadratic differential that is associated to an extremal quasiconformal mapping.

The modern period starts with the works of Gr\"otzsch, Lavrentieff and Ahlfors. We present some of their important ideas in Sections \ref{s:G}, \ref{s:L} and \ref{s:A} respectively.   

Section \ref{s:T} contains information on Teichm\"uller's writings and style. Our aim is to explain why  these writings were not read. Despite their importance, some of his ideas are still very poorly known, even today.

Section \ref{s:T2} is a survey of the main contributions of Teichm\"uller on the use of quasiconformal mappings in the problem of moduli of Riemann surfaces and in problems in function theory.

\medskip

Our digressions also concern the lives and works of several authors we mentioned and who remain rather unknown to quasiconformal theorists. This concerns in particular Nicolas-Auguste Tissot, who had bright ideas that make the link between old problems on conformal and close-to-conformal mappings\footnote{In the present paper, we use the expression ``close-to-conformal" in a broad sense. Generally speaking, we mean by it a mapping between surfaces which is (close to) the best map with respect to a certain desired property. This property may involve angle-preservation (conformality in the usual sense), area-preservation, ratios of distances preservation, etc. Usually the desired property is a combination of these. Making this idea precise is one of the recurrent themes in cartography. Note that we shall also use the notion of ``close-to-conformality"  in the sense of a \emph{conformal map} which minimizes some other distortion parameters, for instance the magnification ratio (see the works of Lagrange, Chebyshev and others that we survey in \S \ref{s:qc}  below).} that arise in cartography, and modern quasiconformal theory.  The  name Tissot is inexistent in the literature on quasiconformal mappings, with the exception of references to his work by Gr\"otzsch in his early papers, see \cite{Groetzsch1930}, and a brief reference to him by Teichm\"uller in his paper \cite{T20}. We also highlight the work of Mikha\"\i l Alekse\"\i evitch Lavrentieff on quasiconformal mappings. This work is well known in Russia, but it is rather poorly quoted in the Western literature on this subject.   
                        
                                           We hope that these historical notes will help the reader  understand the chain of ideas and their development.

There are other historical surveys on quasiconformal mappings; cf. for instance \cite{Lehto1984} and \cite{Caz}.

\part{Geography}

                 \section{Geography and early quasiconformal mappings} \label{s:qc}
     
    We start by recalling the notion of a quasiconformal mapping, since this will be the central element in our paper. 
    
    Let us first consider a linear map $L:V_1\to V_2$ between two 2-dimensional vector spaces $V_1$ and $V_2$ and let us take a circle $C$ in $V_1$, centered at the origin. Its image $L(C)$ is an ellipse in $V_2$, which is a circle if and only if the linear map $L$ is conformal (that is, angle-preserving). We take as a measure of the \emph{distortion} of $L$, or its ``deviation from conformality, " the quantity \[
     D(L)=\log \frac{b}{a},
     \] where $b$ and $a$ are respectively the lengths of the great and small axes of the ellipse $L(C)$.
     The distorsion $D(L)$ does not depend on the choice of the circle $C$.
     
 Let $f:S_1\to S_2$ be now a map between Riemann surfaces, and assume for simplicity that $f$ is of class $C^1$. At each point of $S_1$, we consider the distortion $D(df_x)$ of the  linear map $df_x$ induced by $f$ between the tangent spaces at $x$ and $f(x)$ to $S_1$ and $S_2$ respectively. The \emph{distortion} of $f$ is defined as
     \[D(f)=\sup_{x\in S_{1}} D(df_x)     .\] The value of $D(f)$ might be infinite. The map $f$ is said to be quasiconformal if its distortion is finite. It is conformal if and only if its distortion is zero. 
          
     This definition has a history, and it went through extensive developments. It has also a large number of applications in various fields (complex analysis, differential equations, metric geometry, physics, geography, etc.). This is the subject of the present paper. Our overall presentation is chronological, reflecting the fact that in mathematics the new theories and new results are built on older ones, and the solution of a problem naturally gives rise to new ones.

\medskip
  
  Teichm\"uller explicitly states in his paper \cite{T20} (see \S 164, titled: \emph{Why do we study quasiconformal mappings?}) that the notion of quasiconformal mapping first appeared as a generalization of that of conformal mapping in cartography, and he mentions in this respect the name of  Nicolas-Auguste Tissot (1824--1897). The latter was a French cartographer and mathematics professor at Lyc\'ee Saint-Louis in Paris,\footnote{We note, for the reader who is not familiar with the French educational system, that it was very important to have very good teachers (sometimes they were prominent mathematicians) at the French good lycées, the reason being that these lycées prepared the pupils to the competitive entrance examination at \'Ecole Normale Supérieure or \'Ecole Polytechnique.} and he was also an examiner for the entrance exam of the \'Ecole Polytechnique. The reader is referred to the paper \cite{2016-Tissot} for more information on his life.  Tissot developed a theory in which the distortions of the various projections of the sphere onto the plane that were used in drawing geographical maps play a central role, and he invented a device to measure these distortions. We shall say more about this interesting work, but we start our survey with an exposition of the beginning of conformal geometry and cartography,  and the relation between these two topics. We hope to convey the idea that the idea of quasiconformal mapping is inherent in these works.
 
 Cartography is a science which makes a link between geometry and the real world: ``nature," or ``physis." The Greek word ``geo-graphia"  means "drawing the Earth." Indeed, the classical problems of this field concern the representation on a flat surface of portions of a spherical (or spheroidal) surface. At a practical level, the sphere considered was either the Earth (which was known, since the seventeenth century, to be spheroidal rather than spherical)\footnote{It became generally accepted, starting from the eighteenth century, that the Earth is spheroidal and not spherical, namely, it is slightly flat at the poles. We recall in this respect that the description of the shape of the Earth was a major controversial issue in the seventeenth and eighteenth centuries. This question opposed the English scientists, whose main representative was Newton, to the French, represented by the astronomer Jacques Cassini (1677--1756), the physicist Jean-Jacques d'Ortous de Mairan (1678--1771)  and others, who pretended on the contrary that the Earth is stretched out at the poles.  Newton had concluded in his \emph{Principia} that the Earth is flat at the poles, due to its rotation. More precisely, he expected the flattening to be of the order of 1/230.   Huygens was on the side of Newton.  The controversy on the shape of the Earth is one of the most fascinating scientific issues raised in the eighteenth century science.  It led to long expeditions by French and other scientists to the Lapland and Peru, in order to make precise measurements of the meridians near the poles.} or the heavenly sphere with its multitude of celestial bodies: planets, stars, constellations, etc. A geographical map is a representation in the Euclidean plane of a subset of one of these curved surfaces.\footnote{Of course, at the local level, the Earth cannot be neither spherical nor spheroidal; it has mountains, valleys, canyons, plateaux, etc. In principle, all this has to be taken into account in map drawing. But for the mathematics that is discussed here, it is sufficient to assume that the Earth is spherical. From a practical point of view, this would be perfectly correct if we assume that the map is intended for airplanes that travel a few thousand miles above the surface of the Earth in such a way that they do not face any obstacle.}

\begin{figure}[htbp]
\centering
\includegraphics[width=6cm]{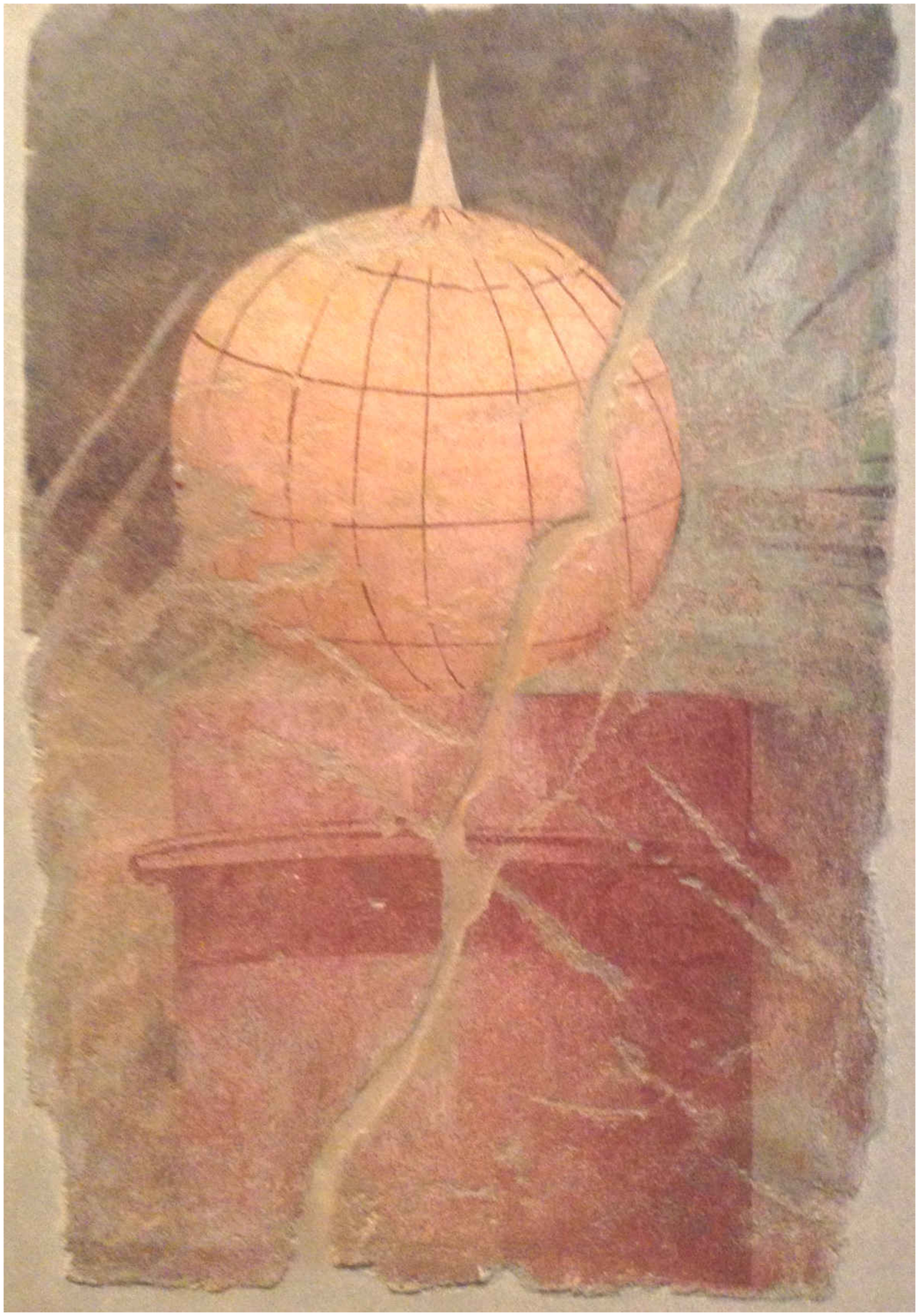}
\caption{\smaller{Celestial globe with coordinates, Roman, c. 50--40 B.C., Metropolitan Museum of Art, New York \emph{(Photo A. Papadopoulos)}.}} \label{Met}
\end{figure}

 Some of the practical questions that arise in cartography led to major problems in abstract mathematics; we shall see this in the next few pages. Furthermore, there is an aesthetic side in the art of map drawing -- this is especially visible in ancient maps -- which certainly complements the beauty of pure geometrical thought. An old representation of coordinates on the Earth by meridians and parallels is reproduced in Figure \ref{Met}.

As a science, the drawing of geographical maps dates back to Greek antiquity. It was known to the mathematicians and geographers of that epoch that there is no representation of the sphere on the plane that is faithful in the sense that proportions of distances are respected.  Despite the fact that no mathematical proof of this statement existed, the idea was widely accepted.  The notion of conformal (angle-preserving) mapping existed. In particular, it was known that areas were always distorted by such a mapping. The relative magnitudes of the various lands needed to be represented in a faithful manner. This came naturally as a need for practical problems linked to harvesting, the distribution of water, calculating taxes, etc. The question of finding  the ``best" maps, that is, those with a minimal amount of ``distortion," with this word referring to some specific property or a compromise between a few desirable properties, arose very naturally. Angle-preservation, or closeness to angle-preservation, was another a desirable property.  Thus, the relation between geography, the theory of quasiconformal mappings, and what we call close-to-conformal mappings has many aspects.
  
  Several prominent mathematicians were at the same time geographers and astronomers. They were naturally  interested in the drawing of geographical maps and they were naturally led to problems related to conformality and close-to-conformality. Among these mathematicians, in Greek Antiquity, the name of Claudius Ptolemy (second century A.D.) comes to the forefront. His work on astronomy and geography is well known, especially through his  \emph{Mathematical syntaxis} and his \emph{Geography}.\footnote{The two majors works of Ptolemy are his \emph{Mathematical syntaxis}, known as the \emph{Almagest}, a comprehensive treatise in thirteen books on astronomy (which also contains chapters on Euclidean and spherical geometry), and his \emph{Geography}, an exposition of all what was known in geography  in the second-century Greco-Roman world. It includes a compilation of geographical maps. The name \emph{Almagest} derives from the Arabic ``al-Majis\b{t}\=\i ,"  which is a deformation of the Greek superlative ``Megistos," meaning ``the Greatest" and referring to the expression ``the Greatest Treatise." The Arabs had the greatest admiration for this work, and several translations, sponsored by the Arabic rulers, were performed starting from the ninth century. A Latin translation was done in the twelfth century, from the Arabic. The \emph{Geography} was translated into Arabic in the ninth century and, later on, from this language into Latin. One of the oldest such translations was made by Herman of Carinthia  in 1143. We refer the interested reader to the discussion and references in the recent paper \cite{Abgrall}.} Besides designing his famous  maps of the Celestial sphere, Ptolemy conceived maps of the Earth, using the mathematical tools he developed in his astronomical works. In his \emph{Geography} (cf. \cite{Ptolemy-geo} and \cite{Ptolemy-Berggren}) he uses the comparison of arcs of circles on  the Celestial sphere with corresponding arcs of circles on the Earth, in order to compute distances. On should mention also Ptolemy's \emph{Planisphaerum}, a treatise  which reached us only through Arabic translations.
    
  Like several mathematicians who preceded him, Ptolemy knew that in drawing  geographical maps, it is not possible to preserve at the same time angles and ratios of distances,
  and he searched for the best reasonable compromise between the corresponding distortions. 
    
 It may be useful to recall in this respect that spherical geometry was developed several centuries before Ptolemy by Greek geometers. One may mention the works of Hipparchus of Nicaea (2nd c. B.C.),  Theodosius of Tripolis (2nd c. B.C.) and Menelaus of Alexandria (1st-2nd c. A.D.). The notion of spherical angle was present in these works. In particular, conformality in the sense of angle-preservation was a meaningful notion.
 
  Ptolemy's \emph{Geography} is divided into eight books. It starts with the distinction between what he calls ``regional cartography" and ``world cartography." Different techniques are used in each subfield, and the goals in each case is different. In regional cartography, there is an aesthetic dominating element, whereas world cartography  involves precise measurements and is used for practical needs. Let us quote Ptolemy (\cite{Ptolemy-Berggren}, Book 1, p. 58):\footnote{The English translation is from the edition of Berggren and Jones \cite{Ptolemy-Berggren}.}
  \begin{quote}\small
  The goal of regional cartography is an impression of a part, as when one makes an image of just an ear or an eye; but [the goal] of world cartography is a general view, analogous to making a portrait of the whole head. [...]
  
  Regional cartography deals above all with the qualities rather than the quantities of the things that it sets down; it attends everywhere to likeness, and not so much to proportional placements. World cartography, on the other hand, [deals] with the quantities more than the qualities, since it gives consideration to the proportionality of distances for all things, but to likeness only as far as the coarser outlines [of the features], and only with respect to mere shape. Consequently, regional cartography requires landscape drawing, and no one but a man skilled in drawing would do regional cartography. But world cartography does not [require this] at all, since it enables one to show the positions and general configuration [of features] purely by means of lines and labels.
  
  For these reasons, [regional cartography] has no need of mathematical method, but here [in world cartography] this element takes absolute precedence. Thus the first thing that one has to investigate is the Earth's shape, size, and position with respect to its surroundings [i.e., the heavens], so that it will be possible to speak of its own part, how large it is and what it is like, and moreover [so that it will be possible to specify] under which parallels of the celestial sphere each of the localities of this [known part] lies. From this last, one can also determine the lengths of nights and days, which stars reach the zenith or are always borne above or below the horizon, and all the things that we associate with the subject of habitations. 
  
  These things belong to the loftiest and loveliest of intellectual pursuits, namely to exhibit to human understanding through mathematics [both] the heavens themselves in their physical nature (since they can be seen in their revolution about us), and [the nature of] the Earth through a portrait (since the real [Earth], being enormous and not surrounding us, cannot be inspected by any one person either as a whole or part by part).
    \end{quote}
   Thus, regional cartography has a taste of topology. Distortions, whether at the level of angles, distances or areas, play a negligible role in this art. For a modern example of a map falling into Ptolemy's category of regional cartography,  one may think of a railroad of subway network map, where the information that is conveyed is only contained in the lines and their intersections. In fact, in such a map, it is the topology of the network that is the interesting information.
  
  Ptolemy's writings contain valuable information on works done by his predecessors. It is important to remember that Ptolemy, like Anaximander (c. 610 --546 B.C.), Eratosthenes (3d c. B.C.),\footnote{Eratosthenes  gave a measure of the circumference of the Earth which is correct up to an error of the order of 1\%.} Hipparchus of Nicaea (c. 190--c. 120 B.C.), Marinus of Tyre (c. 70--130 A.D.) and the other Greek mathematicians and geographers that he quotes, considered that the Earth has a spherical shape. The myth of a flat Earth had very little supporters among scientists in Greek antiquity. 
In his \emph{De Caelo} \cite{Aristotle-Heavens} (294 b14-15), Aristotle reports that the pre-Socratic philosophers Anaximenes (6th. century B.C.), who highlighted the concept of infinite, his student Anaxagoras (5th. century B.C.), who was the first to establish a philosophical school in Athens, and Democritus (4th--5th. century B.C.), who formulated an atomic theory of the universe, already considered that the Earth is spherical. Plato added  philosophical reasons to the belief in a spherical Earth, namely, he considered that the Demiurge could only give to the material world the most symmetrical form, that is, the spherical one. In the \emph{Timaeus} (33b) \cite{Plato}, we read:
  \begin{quote}\small
  He wrought it into a round, in the shape of a sphere, equidistant in all directions from the center to the extremities, which of all shapes is the most perfect and the most self-similar, since he deemed that the similar is infinitely fairer than the dissimilar. And on the outside round about, it was all made smooth with great exactness, and that for many reasons. 
  \end{quote}
  
  Let us note by the way that the eleventh-century Arabic polymath Ibn al-Haytham, at the beginning of his treatise on  isoepiphany, declares that the reason for which the entire universe and the Earth have a spherical shape is that the sphere has the largest volume among the solid figures having the same area. (See the translation of Ibn al-Haytham's treatise in \cite{Rashed1}, English version,  p. 305--342.)

  The stereographic projection, that is, the radial projection of the sphere from a point on this surface (say the North pole) onto a plane passing through the center and perpendicular to the radius passing through this point (in the case considered, this would be the equator plane), which has probably always been the most popular projection of the sphere among mathematicians, was already used by Hipparchus back in the second century B.C.; cf. Delambre \cite{Delambre}, Vol. 1, p. 184ff.,   d'Avezac \cite{Avezac} p. 465 and Neugebauer \cite{Neuge-astro} p. 246. Let us also add that to Hipparchus is also attributed the invention of the astrolabe, a multipurpose hand-held geographical instrument which made out of several metallic plates that fit and roll into each other and where each of them is considered as a planar representation of the celestial sphere.  The astrolabe allows its user to  locate himself (altitude and longitude), know the hour of the day or of the night and find the location of the main stars and planets. 
    The instrument was used by astronomers and navigators. It became very popular in the Islamic world, starting from the eighth century,  because it allowed its user to find, at a given place, the exact lunar time and the direction of the Mecca for everyday prayer. Obviously, nontrivial geometrical problems arise in the conception of the astrolabe, since this requires the  transformation of spherical periodic movements of the celestial sphere (the global movement of the whole  sphere from East to West around an axis which is perpendicular to the equator, the movement of the sun and of the planets on the ecliptic circle, etc.) into plane geometrical ones.  Since Greek antiquity, the mathematical theory of the astrolabe involved several aspects of geometry: Euclid's plane and solid geometry, the theory of proportions, Apollonius' theory of conics, and descriptive geometry. The use of the astrolabe disappeared with the popularization of modern clocks. There are basically three kinds of writings on the astrolabe:  its confection, its use, and the mathematical theory that is at its basis. Sereval Arabic scholars, including al-\d{S}\=agh\=an\=\i , al-Q\=uh\=\i , Ibn Sahl and Ibn `Ir\=aq, described the astrolabe in several writings, and they transmitted it to the West.  The interested reader may consult the book \cite{Rashed-Geometry} by Rashed.  
  
  The stereographic projection was not the most useful one for the drawing of geographical maps, because the distortion of distances and areas becomes very large at points which are far from the equator.  Better projections were known to geographers; we shall mention some of them in what follows. Neugebauer, in his survey \cite{Neu1948} (p. 1037--1039) mentions a cylindrical projection  by Marinus of Tyre in which distances are preserved along all the meridians and along the parallel passing by the island of Rhodes (36${}^{\mathrm{o}}$).\footnote{This is the parallel which was situated at the middle of the civilized world, from Spain to China.}  In this projection,  distances along the parallels which are north of Rhodes are expanded, whereas those along the parallels south of this island are contracted. Marinus' maps are considered as models for the Mercator maps, which appeared fourteen centuries later. It is also known that Christopher Columbus, in the preparation of his journeys, relied in part on Marinus' maps and computations, in part contained in Ptolemy's \emph{Almagest}. Ptolemy introduced conical projections which in some cases are improvements of the one of Marinus. In these projections, distances along all meridians are preserved, and they are also preserved along the parallels which pass through Rhodes and through two other meridians: Thule (63${}^{\mathrm{o}}$) and the equator (0${}^{\mathrm{o}}$).
  We also mention that Marinus of Tyre, besides being a  geographer, was, like Ptolemy, a mathematician. In Chapters XI to XX of his \emph{Geography}  (cf. \cite{Ptolemy-geo} p. 40ff.), Ptolemy criticizes Marinus'  cylindrical projections because of their significant amount of  distortion of lengths.

During the middle ages, science was flourishing in the Arabic world, and geography was part of it. The Arabs studied the works of Marinus of Tyre and Ptolemy, and they built upon them a whole school of geographical map drawing, including valuable new developments. One should mention at least the works of the historian, traveler and geographer al-Mas`\=ud\=\i \ (c. 896--956), and those of the mathematician, astronomer, philosopher, geographer, physician and pharmacologist al-B\=\i r\=un\=\i  \ (973 --1048) who wrote an important  treatise in which he describes several new projections; cf. \cite{Richter}.

\begin{figure}[htbp]

\centering
\includegraphics[width=12cm]{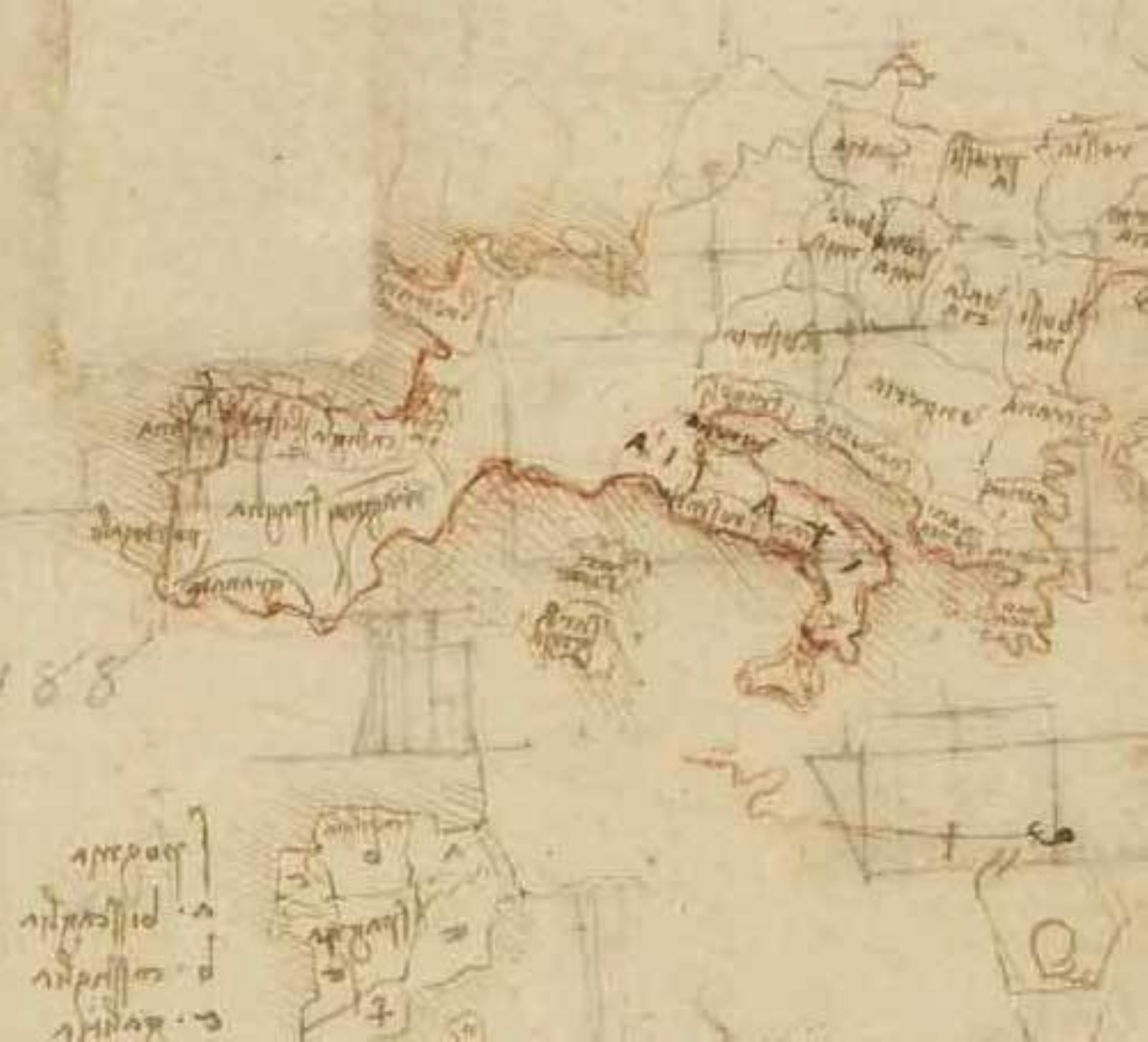}
\caption{\smaller{A map drawing by Leonardo da Vinci.}} \label{Leonardo}
\end{figure}

   At the Renaissance, the need for drawing geographical maps was strongly motivated by the discovery of new lands. Leonardo da Vinci and Albrecht D\"urer, who, besides being artists, were remarkable mathematicians, were highly interested in geography. They worked on the technical problems raised by the drawing of geographical maps, cf. \cite{Oberhummer}, \cite{Leonardo} and \cite{Panofsky}, and they also valued the beauty of these drawings. See Figure \ref{Leonardo}.

   In the modern period, the works of Euler and Gauss on geography are well known, but there are also many others.  Mathematical questions related to conformality and close-to-conformality that arise from the drawing of geographical maps or that are directly related to that art were thoroughly studied by Lambert, Lagrange, Chebyshev and Beltrami. Furthermore, in the works of several French geometers including Liouville, Bonnet, Darboux and others, cartography became an application of infinitesimal calculus and an important chapter of the differential geometry of surfaces. We shall briefly review some of these important works. They constitute a significant part of the history of quasiconformal mappings.

  We start with Euler, then talk about  Lambert, Lagrange, Chebyshev, Gauss, Liouville, Bonnet, Darboux and Beltrami, before reviewing the pioneering work of Tissot.
  
  Euler, at the Academy of Sciences of Saint Petersburg, besides being a mathematician, had the official position of cartographer. He participated to the huge project of drawing geographical maps of the new Russian Empire. Having  reliable maps was important to the new Russian rulers, and they were used for several purposes: determining the precise boundary of the empire, determining areas of pieces of lands for taxation, irrigation, and other purposes, giving precise distances between cities, with the lengths of the various roads joining them, etc.
It is interesting to know that Euler was so much involved in this activity that he considered (or at least he claimed) that it was responsible for the deterioration of his vision. In a letter to Christian Goldbach, dated August 21st (September 1st, new Style), 1740, he writes (cf.  \cite{Euler-Goldbach} p. 163, English translation p. 672):
\begin{quote}\small
Geography is fatal to me. As you know, Sir, I have lost an eye working on it; and just now I nearly risked the same thing again. This morning I was sent a lot of maps to examine, and at once I felt the repeated attacks. For as this work constrains one to survey a large area at the same time, it affects the eyesight much more violently than simple reading or writing. I therefore most humbly request you, Sir, to be so good as to persuade the President by a forceful intervention that I should be graciously exempted from this chore, which not only keeps me from my ordinary tasks, but also may easily disable me once and for all.  I am with the utmost consideration and most respectfully, Sir, your most obedient servant.
 \end{quote}

Among Euler's works on geography, we mention his \emph{Methodus viri celeberrimi Leonhardi Euleri determinandi gradus meridiani pariter ac paralleli telluris, secundum mensuram a celeb. de Maupertuis cum sociis institutam}
(Method of the celebrated Leonhard Euler for determining a degree of the meridian, as well as of a parallel of the Earth, based on the measurement undertaken by the celebrated de  Maupertuis and his colleagues) \cite{E132}. This work was presented to the Academy of Sciences of Saint-Petersburg in 1741 and published in 1750.\footnote{There was sometimes a significant lapse of time between the time when Euler wrote a paper and the time when the paper appeared in print. One reason for this delay is that the journals of the Academies of Saint Petersburg  and of Berlin, where most of Euler's papers were published, received a large number of such papers, and, since the number of pages of each volume was limited, the backlog became gradually substantial. In fact, several papers of Euler were presented to these Academies even after his death and  published much later. Starting from the year 1729, and until 50 years after Euler's death, his works continued to fill about half of the scientific part of the \emph{Actes} of the Saint-Petersburg Academy. Likewise, between the years 1746 and 1771, almost half of the scientific articles of the \emph{M\'emoires}  of the Academy of Berlin were written by Euler.
It is also known that in some cases, Euler purposely delayed  the publication of his memoirs, in order to leave the primacy of the discoveries to a student or young collaborator working on the same subject, especially when he considered that their approach was at least as good as his. A famous instance  of this generous attitude is when Euler  waited so that  Lagrange, who was nineteen years younger than him, finishes and publishes his work on the calculus of variations, before he sends his own work to his publisher. This was pointed out by several authors. The reader can refer for instance to Euler's obituary by Condorcet \cite{Condorcet-Eloge} p. 307.} 

From a more theoretical point of view, Euler published in 1777 three memoirs on mappings from the sphere onto the Euclidean plane:
\begin{enumerate}
\item  \emph{De repraesentatione superficiei sphaericae super plano} (On the representation of  spherical surfaces on a plane) \cite{Euler-rep-1777};
\item
  \emph{De proiectione geographica superficiei sphaericae} (On the geographical projections of spherical surfaces) \cite{Euler-pro-1777};
\item  \emph{De proiectione geographica Deslisliana in mappa generali imperii russici usitata} (On Delisle's geographic projection used in the general map of the Russian empire)  \cite{Euler-pro-Desli-1777}. 
\end{enumerate}
These three memoirs were motivated by the practical question of drawing geographical maps, and they led to mathematical investigations on conformal and close-to-conformal projections. The title of the third memoir  refers to Joseph-Nicolas Delisle (1688-1768), a leading French geographer, astronomer, meteorologist, astrophysicist, geodesist, historian of science and orientalist who worked, together with Euler, at the Saint Petersburg Academy of Sciences.\footnote{Delisle is one of the eminent scientists from Western Europe who were invited by Peter the Great at the foundation of the  Russian Academy of Sciences. The monarch  signed the foundational decree of his Academy on February 2, 1724. The group of scientists that were present at the opening ceremony included the mathematicians Nicolaus and Daniel Bernoulli and Christian Goldbach. Delisle joined the Academy in 1726, that is, two years after its foundation, and one year before Euler.  Delisle was in charge of the observatory of Saint Petersburg, situated on the Vasilyevsky Island. This observatory was one of the finest in Europe. Thanks to his position, Delisle had access to the most modern astronomical instruments. He was in charge of drawing maps of the Russian empire. Euler assisted him at the observatory and in his map drawing. Delisle stayed in  Saint Petersburg from 1726 to 1747. In the first years following his arrival to Russia, Euler helped Delisle in recording astronomical observations which were used in meridian tables. Delisle returned to Paris in 1747, where he founded the famous observatory at the H\^otel de Cluny, thanks to a large amount of money he had accumulated in Russia. He also published there a certain number of papers containing information he gathered during his various journeys inside Russia and in the neighboring territories (China and Japan). It seems that some of the information he published was considered as sensible, and because of that, Delisle acquired the reputation of having been a spy.} We now review important ideas contained in these three memoirs.
 
 The first memoir \cite{Euler-rep-1777} contains important results and techniques that include the questions on geographical maps in the setting of differential calculus and the calculus of variations. In the introduction, Euler declares that he does not consider only maps obtained by central projection of the sphere onto the plane which are subject to the rules of perspective,  (he calls the latter ``optical projections"), but he considers the mappings ``in the widest sense of the word."  Let us quote him (\S 1 of \cite{Euler-rep-1777}):
 \begin{quote}\small
 In the following I consider not only optical projections, in which the different points of a sphere are represented on a plane as they appear to an observer at a specified place; in other words, representations by which a single point visible to the observer is projected onto the plane by the rules of perspective; indeed, I take the word ``mapping" in the widest possible sense; any point of the spherical surface is represented on the plane by any desired rule, so that every point of the sphere corresponds to a specified point in the plane, and inversely; if there is no such correspondence, then the image of a point from the sphere is imaginary.
 \end{quote}

 In \S\,9, Euler proves that there is no ``perfect" or ``exact" mapping from the sphere onto a plane, that is, a mapping which is distance-preserving (up to some factor) on some set of curves. This result is obtained through a study of partial differential equations. The precise statement of Euler's theorem, with a modern proof based on Euler's ideas, are given in the paper \cite{CP}.

 With this negative result established, Euler says that one has to look for best approximations.  He writes:  ``We are led to consider representations which are not similar, so that the spherical figure differs in some manner from its image in the plane." He then examines several particular projections of the sphere, searching systematically searches for the partial differential equations that they satisfy. 
In doing so, he highlights the following three kinds of maps:

\begin{enumerate}
\item \label{map1} Maps where the images of all the meridians are perpendicular to a given axis (the ``horizontal" axis in the plane), while all parallels are sent parallel to it. 

\item \label{map3} Maps which ``preserve the properties of small figures,"  which, in his language, means conformal.

\item \label{map2} Maps where surface area is represented at its true size.

\end{enumerate}
Euler gives examples of maps satisfying each of the above three properties. These maps are obtained by various ways: projecting the sphere onto a tangent plane, onto a cylinder  tangent to the equator, etc. Euler then studies distance and angle distortion under these various maps. At the end of his memoir (\S 60), he notes that his work has no immediate practical use:\footnote{For the three papers of Euler on geography, we are using the translation by George Heine.}
\begin{quote}\small
 In these three hypotheses is contained everything ordinarily desired from
geographic  as  well  as  hydrographic  maps.   The  second  hypothesis  treated
above even covers all possible representations.  But on account of the great
generality of the resulting formulae,  it is not easy to elicit from them any
methods of practical use.  Nor, indeed, was the intention of the present work
to go into practical uses, especially since, with the usual projections, these
matters have been explained in detail by others.
\end{quote}
 
 In contrast, in the memoir \cite{Euler-pro-1777}, Euler studies projections that are useful for practical applications. He  writes (\S 20):
 \begin{quote}\small
 Moreover, let it be remarked, that this method of projection is ``extraordinarily appropriate" for the practical applications required by Geography,
for it does not distort too much any region of the Earth. It is also important
to note that with this projection, not only are all meridians and circles of
parallels exhibited as circles or as straight lines, but all great circles on the
sphere are expressed as circular arcs or straight lines. 
\end{quote}

The memoir \cite{Euler-pro-Desli-1777} also contains important ideas. Euler starts by reviewing the main properties of the stereographic projection: the images of the parallels and meridians intersect at right angles, and it is conformal (he says, a similitude on the small scale). He then exhibits the inconveniences of this projection: length is highly distorted on the large, especially if one has to draw maps of large regions of the Earth. He gives the example of the map of Russia in which the province of Kamchatka is distorted by a factor of four, compared to a region on the center of the map.  Another disadvantage he mentions is the distortion of the \emph{curvature} of the meridians, even if they are sent to circles, and he again gives the example of the Kamchatka meridian.  He then considers another possible projections which is commonly used, where the meridians are sent to straight lines which intersect at the image of the pole, and which presents similar disadvantages. Euler then explains the advantage of a new map elaborated by Delisle, and he develops a mathematical theory of this map. 
It is interesting to recall how Euler formulates the problem, as a problem of ``minimizing the maximal error" over an entire region:
\begin{quote}\small
Delisle, the most celebrated astronomer and geographer of the time, to whom the care
of such a map was first entrusted, in trying to fulfill these conditions, made the relationship between latitude and longitude exact at two noteworthy parallels. He was of the opinion that if the named circles of parallel were at the same distance from the middle parallel of the map as from its outermost edges, the deviation could nowhere be significant. Now the question is asked, which two circles of parallel ought to be chosen, so that the maximum error over the entire map be minimized.
\end{quote}

In developing the mathematical theory of Delisle's map, Euler finds that one of their advantages is that while meridians are represented by straight lines, the images of the other great circles ``do not deviate considerably from straight lines" (\S 22 of Euler's paper \cite{Euler-pro-Desli-1777}). Such a condition will appear several decades later (almost a century), in a paper by Beltrami \cite{Beltrami1865}, which we review below. It is probably the first time where we encounter,  formulated in a mathematical setting, an early notion of ``quasigeodesic." Euler then investigates the way the image of a great circle on the map differs from a straight line (\S 23ff.)
The conclusion of the memoir is:
\begin{quote} \small
 In this projection is obtained the extraordinary advantage, that straight lines, which go from any point to any other point, correspond rather exactly to great circles and therefore the distances between any places on the map can be measured by using a compass without considerable error. Because of these important characteristics the projection discussed was preferred before all others for a general map of the Russian Empire, even though, under rigorous examination, it differs not a little from the truth.

\end{quote}

Euler had several young collaborators and colleagues  working on cartography. We mention in particular F. T. Schubert,\footnote{Friedrich Theodor von Schubert (1758-1825), like Euler,  was the son of a protestant pastor. His parents, like Euler's parents, wanted their son to study theology.  Instead, Schubert decided to study mathematics and astronomy, and he did it without the help of teachers. He eventually became a specialist in these two fields,  and he taught them as a private teacher during his stays outside his native country (Germany). After several stays in different countries,  Schubert was appointed, in 1785, assistant at the Academy of Sciences of Saint Petersburg, at the class of geography. This was two years after Euler's death. In 1789 he  became full member of that Academy, and in 1803  director of the astronomical observatory of the same Academy.} one of his direct followers who became a specialist of spherical geometry, a field closely related to cartography and on which Euler published several memoirs. We refer the reader to the article \cite{Papa-Inde2} for more information about Schubert. Schubert's papers on cartography  include \cite{S0},  \cite{S00} and \cite{S01}. The forthcoming book \cite{Caddeo2} contains an analysis of the works of Euler and his collaborators on spherical geometry.

  After Euler, we talk about his colleague and fellow countryman J. H. Lambert,\footnote{Depending on the biographers, the Alsatian self-taught mathematician Johann Heinrich Lambert  (1728-1777) is considered as French, Swiss or German. He was born in Mulhouse, today the second largest French city in Alsace, but which at that time was tied to the Swiss Confederation by a treaty which guaranteed to it a certain independence from the Holy Roman Empire and a certain relative peace away from the numerous conflicts which were taking place in the French neighboring regions. 
Lambert belonged to a French protestant family which was exiled after the revocation by Louis XIV of the Edict of Nantes. We recall that this edict of Nantes was issued (on  April 13, 1598) by the king Henri IV of France at the end of the so called Wars of Religion which devastated France during the second half of the 16th century. Its goal was to  grant the Protestants a certain number of rights, with the aim of promoting civil unity. Louis XIV suspended the Edict of Nantes on October 22, 1685, by another edict, called he Edict of Fontainebleau.  At the same time, he ordered the closing of the Protestant churches and the destruction of their schools.  This caused the exile of at least 200 000 protestants out of France.
 
 Euler had a lot of consideration for his compatriot's work and he recommended him for a position at the Berlin Academy of Sciences, where Lambert was hired and where he spent the last ten years of his life. We mention by the way that Lambert was one of the most brilliant precursors of hyperbolic geometry, and in fact, he is considered as the mathematician who came closer to discovering it, before Gauss, Lobachevsky and Bolyai. Indeed, in his {\it Theorie der Parallellinien} (Theory of parallels), written in 1766, Lambert developed the bases of a geometry where all the Euclidean postulates hold except the parallel postulate, and where the latter is replaced by its negation. In developing this geometry, his aim was to find a contradiction, since he believed, like all the other geometers before him, that such a geometry cannot exist. The memoir ends without conclusion, and Lambert did not publish it, the conjectural reason being that Lambert was not sure at the end whether such a geometry exists or not. The memoir was published, after Lambert's death, by Johann III Bernoulli. We refer the reader to \cite{Lambert-Blanchard} for the first translation of Lambert's  {\it Theorie der Parallellinien} together with a mathematical commentary, and to the paper \cite{Papa-Inde3} for an exposition of the main results in his memoir.} who is considered as the founder of modern cartography.  His \emph{Anmerkungen und Zus\"atze zur Entwerfung der Land- und Himmelscharten} (Remarks and complements for the design of terrestrial and celestial maps, 1772) \cite{Lamb-Anmer} contains seven new geographical maps, each one having important features. They include (among others) the so-called Lambert conformal conic projection, the transverse Mercator, the Lambert azimuthal equal area projection, and the Lambert cylindrical equal-area projection. 
The Lambert projections are mentioned in the modern textbooks of cartography, and some of them were still in use, for military and other purposes, until recent times (that is, before the appearance of satellite maps).  For instance, the \emph{Lambert conformal conic projection} was adopted by the French artillery during the First Wold War. It is obtained by projecting conformally  the surface of the Earth on a cone that touches the sphere along a parallel. The images of the meridians are concurrent lines and those of the parallels are circles centered at the intersection points of the meridians. When the cone is unrolled it has the form of a fan. The distorsion is minimized between two parallels within which lies the region of interest. In the same memoir \cite{Lamb-Anmer}, Lambert gave a mathematical characterization of an arbitrary angle-preserving map from the sphere onto the plane. Euler later gave another solution of the same question. Lambert, in his work, takes into account the fact that the Earth is spheroidal and not spherical. The memoir  \cite{Lamb-Anmer} is part of the larger treatise \emph{Beitr\"age zum Gebrauche der Mathematik und deren Anwendung durch} (Contributions to the use of mathematics and its applications) \cite{Lambert-Bey}.
   
  Lagrange, whose name is associated to that of Euler for several mathematical resulrs, published two memoirs on cartography, \cite{Lagrange1779}, two years after those of Euler appeared in print. In these memoirs, Lagrange extends some works by Euler and Lambert. In the introduction of his memoir  \cite{Lagrange1779}, p. 637, he writes:\footnote{This and the other translations from the French in this paper are ours.}
\begin{quote}\small
A geographical map is nothing but a plane figure which represents the surface of the Earth, or some part of it. [...] But the Earth being spherical, or spheroidal, it is impossible to represent on a plane an arbitrary part of its surface without altering the respective positions and distances of the various places; and the greatest perfection of a geographical map will consist in the least alteration of these distances.
\end{quote}
Lagrange then surveys various projections of the Earth and of the Celestial sphere that were in use at his time. Among them is the stereographic projection, with its two main properties, which he recalls:

\begin{enumerate}
\item 
 The images of the circles of the sphere are circles of the plane. Thus, in particular, if we want to find the image of a meridian or a parallel, then the image of three points suffices. 
\item The stereographic projection is angle-preserving, a property which is a consequence of the first one.
\end{enumerate}

Lagrange attributes the stereographic projection to Ptolemy (\cite{Lagrange1779} p. 639) who used it in the construction of astrolabes and celestial planispheres. He says that Ptolemy was aware of the first property, and that he describes it in a treatise known as the \emph{Sphaerae a planetis projectio in planum}.\footnote{This the is famous \emph{Planisphaerum} which we mentioned in \S \ref{s:qc}.} Latin version is translated from the Arabic;  the Greek original is lost.) Lagrange says that the angle-preserving property of the stereographic projection may not have been noticed by the Greek astronomer. He then introduces other conformal projections, and he notes that there are infinitely many such maps. In fact, after introducing the great variety of maps from the sphere onto a plane that are projections from a point (allowing also this point to be at infinity), Lagrange considers more general maps. He declares, like Euler before him, that one may consider geographical maps which are arbitrary representations of the surface of the sphere onto a plane. Thus, he is led to the question:
\begin{quote}
 \emph{What is a general mapping between two surfaces?} 
 \end{quote}
 The question is a two-dimensional analogue of the question \begin{quote}
 \emph{What is a function?}  
 \end{quote}
  which was an important issue in eighteenth-century mathematics. It is related to the question of whether a general function is necessarily expressed by an analytic formula or not, which gave rise to a fierce debate involving the most prominent mathematicians of that time, and which lasted several decades; see e.g. the exposition in \cite{Papa-Physics}.

  Returning to the maps from the sphere onto the Euclidean plane, Lagrange writes  (p. 640) that ``the only thing we have to do is to draw the meridians and the parallels according to a certain rule, and to plot the various places relatively to these lines, as they are on the surface of the Earth with respect to the circles of longitude and latitude." Thus, the images of the meridians and the parallels are no more restricted to be circles or lines. They can be, using Lagrange's terms,  ``mechanical lines," that is, lines drawn by any mechanical device. We note that this notion of curve as a mechanical line exists since Ancient Greece.  (In fact, this is the way curves were defined.) It is known today that these are the most general curves that one can consider; cf.  the result of Kempe  that any bounded piece of an algebraic curve is drawable by some linkage \cite{Kempe} and its wide generalization by W. P. Thurston, and by Nikolay Mnev who proved a conjecture of Thurston on this matter \cite{Mnev}. This work is reported on by Sossinsky in \cite{Sossinsky}. Regarding the generality of the mappings from the sphere onto the plane, Lagrange refers to the work of Lambert, who was the first to study arbitrary angle-preserving maps from the sphere to the plane, and who, in his memoir \cite{Lambert-Bey}, solved the problem of characterizing a least-distortion map among those which are angle-preserving. Lagrange declares that Lambert was the first to consider arbitrary maps from the sphere to the plane. He recalls that Euler, after Lambert, gave a solution of the problem of minimizing distortion of an arbitrary angle-preserving map, and he then gives his own solution, by a method which is different from those of Lambert and Euler. He considers in detail the case where the images of the meridians and the parallels are circles. He solves the problem of finding all the orthonormal projections of a surface of revolution which send meridians and parallels to straight lines or circles.

We highlight a formula by Lagrange that gives the local distortion factor of  a conformal map. This quantity, which Lagrange denotes by $m$, is a function on the spherical (or the spheroidal) surface which measures the local distortion of length at a point. It is defined as the ratio of the infinitesimal length element at the image by the infinitesimal length element at the source. The fact that the map is conformal map makes this quantity well defined.
Is is given by the following formula (using Lagrange's notation, p. 646 of Vol. IV of his \emph{Collected works} \cite{Lagrange1779}):
\begin{equation}\label{distortion}
m=\frac{\sqrt{f'(u+t\sqrt{-1})F'(u-t\sqrt{-1})}}{\sin s},
\end{equation}
where, in the case where the Earth is considered as spherical, $(s,t)$ are the coordinates such that 

-- $s$ is the length of the arc of meridian counted from the pole;

-- $t$ is the angle which the meridian makes with a fixed meridian;

-- $u= \log \big( k\tan \frac{s}{2}\big)$;

-- $f$ and $F$ are arbitrary functions.

Lagrange then gives a formula for the distortion of maps which send meridians and parallels to circles or circles and straight lines. 
These formulae were used later on by Chebyshev in his work on geographical maps which we review below. 
 We shall return below to Lagrange's work, when we shall talk about Chebyshev. In the meanwhile, let us talk about Gauss.
  
   Gauss, in the preface of his paper \cite{Gauss-Copenhagen}, declares that his aim is only to construct geographical maps and to study the general principles of geodesy for the task of land surveying. Indeed, it was this task that led him gradually to the investigation of surfaces and their triangulations,\footnote{Triangulations of surfaces originate in cartography. The idea is that to know the position of point, it is often practical to measure the angles in a triangle whose vertices are this point and two other reference points, and use (spherical) trigonometric formulae. In such a triangle, two angles and the length of a side contained by them are known, therefore the three side lengths are completely determined. We quote Hewitt from \cite{Hewitt} (Chapter 3): ``Triangulation had first emerged as a map-making method in the mid sixteenth century when the Flemish mathematician Gemma Frisius set out the idea in his \emph{Libellus de locorum describendorum ratione} (Booklet concerning a way of describing places), and by the turn of the eighteenth century it had become the most respected surveying technique in use. It was a similar method to plane-table surveying, but its instruments conducted measurements over much longer distances."}  to the method of least squares (1821), and then to his  famous treatise \emph{Disquisitiones generales circa superficies curvas} (General investigations on curved surfaces) (1827). In the latter, Gauss gives examples extracted from his own measurements. He writes, in \S 27 (p. 43 of the English translation \cite{Gauss-English}): ``Thus, e.g., in the greatest of the triangles which we have measured in recent years, namely that between the points Hohehagen, Brocken, Inselberg, where the excess of the sum of the angles was 14."85348, the calculation gave the following reductions to be applied to angles: 
Hohehagen: 4."95113;
Brocken: 4."95104;
Inselberg: 4"95131." In the same work, Gauss proves the following, which he calls ``remarkable theorem" (\emph{Theorema Egregium}) that explains in particular why curvature is the obstruction for a sphere to be faithfully represented on the plane (\S 12; p. 20 of the English translation)
\begin{quote} \small
\emph{If the curved surface is developed upon any other surface whatever, the measure of curvature in each point remains unchanged.}
\end{quote}

    In  1825, Gauss published a paper in the \emph{Astronomische Abhandlungen}\footnote{This paper won a prize for a competition proposed by the Copenhagen Royal Society of Sciences. The subject of the competition was:
 ``To represent the parts of a given surface onto another surface in such a way that the representation is similar to the original in its infinitesimal parts." 
  The prize was set by Heinrich Christian Schumacher, a famous German-Danish astronomer in Copenhagen, who was a friend of Gauss and who had been his student in G\"ottingen. A letter from Gauss to Schumacher dated July 5, 1816 shows that the solution to this question was already known to Gauss at that time; cf. Gauss's \emph{Werke} Vol. 8, p. 371. In fact, it was Gauss who, in that letter, proposed to Schumacher to set this problem as the prize contest.  Schumacher responded positively, and the question became the subject of the Academy contest. It is generally assumed that Schumacher had no doubt that the prize will go to Gauss. An English translation of Gauss's paper was published in the \emph{Philosophical magazine} (1828), with the subtitle: \emph{Answer to the prize question proposed by the Royal Society of Sciences at Copenhagen}. It is reproduced in \cite{Smith}, Volume 3. There is also a French translation \cite{Gauss-s}.} (Memoirs on astronomy) whose title is \emph{Allgemeine Aufl\"osung der Aufgabe, die Teile einer gegebenen Fl\"ache auf einer andern gegebenen Fl\"ache so abzubilden dass die Abbildung dem Abgebildeten in den kleinisten Theilen \"ahnlich wird.} (General solution of the problem: to represent the parts of a given surface on another so that  the smallest parts of the representation shall be similar to the corresponding parts of the surface represented) \cite{Gauss-Copenhagen}. In this paper, Gauss shows that every sufficiently small neighborhood  of a point in an arbitrary real-analytic surface can be mapped conformally onto a subset of the plane.  

We mention two other important papers by Gauss related to geodesy:  \emph{Bestimmung des Breitenunterschiedes zwischen den Sternwarten von G\"ottingen und Altona durch Beobachtungen am Ramsdenschen Zenithsector} (Determination of the latitudinal difference between
the observatories in G\"ottingen and Altona by observations with a
Ramsden zenith sector) \cite{Gauss-Bestimmung} (1928) and \emph{Untersuchungen \"uber Gegenst\"ande der h\"ohern Geod\"asie} (Research on objects of higher geodesy) \cite{Gauss-Untersuchungen} (1843 and 1847). It is in the latter that Gauss uses the method of least squares.

One may also talk about Riemann, who was a student of Gauss. His famous Riemann Mapping Theorem, and its later culmination, the Uniformization Theorem, concern the existence of angle-preserving maps between simply connected  Riemann surfaces. Reviewing these developments would be another story, on which several good survey articles exist. We refer the reader to the book \cite{S-Gervais}.

The works of Euler, Lagrange and Gauss are reviewed in the second doctoral thesis\footnote{In France, a doctoral thesis had two parts, a \emph{first thesis} and a \emph{second thesis}.  The subject of the second thesis was proposed by the jury, a few weeks before the thesis defense. It was not necessary that this work be original; generally its was an exposition of some topic (different from  the first thesis's subject). The second thesis disappeared with the reform of the doctoral studies in France, in the 1990s.} of Ossian Bonnet \cite{Bonnet-these}, defended in 1852. Bonnet declares in the introduction that the aim of his thesis is to simplify the works of the three authors on the problem, which he says was first addressed by Lambert, of  characterizing all maps between surfaces which are ``similarities at the infinitesimal level." Thus, we are led once again to angle-preserving maps. Euler and Lagrange, he says, gave a solution of the problem of characterizing these maps, under the condition that the surface is a sphere or, at most, a surface of revolution, and Gauss gave the solution, in his Copenhagen-prize memoir,  for maps between arbitrary surfaces.

The introduction of the thesis contains a short historical survey of cartography. The techniques that Bonnet uses are those of the differential geometry of surfaces. The thesis was published in the \emph{Journal de math\'ematiques pures et appliqu\'ees}. At the end of the paper, Liouville, who was the editor of the journal, declares in a note that the questions addressed in this paper were treated in the lectures he gave at the \emph{Coll\`ege de France} during the academic year 1850--1851, and that he hopes to publish these lectures some day. He also recalls that he explained part of his ideas on the subject in the Notes of his edition of Monge's \emph{Application de l'Analyse \`a la G\'eom\'etrie} \cite{Monge}. These notes, all written by Liouville,  contain the solution of the question he calls (Note V) the \emph{three-dimensional geographical drawing problem}. The title is \emph{Du trac\'e g\'eographique des surfaces les unes sur les autres} (On the geographical drawing of surfaces one onto the other). In the introduction of that note, Liouville formulates again the problem as the one of finding a mapping between two surfaces which is a similarity at the infinitesimal level. He says that this amounts to requiring that each infinitesimal triangle on the first surface is sent by the map to a similar infinitesimal triangle on the second one. He then formulates the same problem 
using the infinitesimal length elements: the ratio between each infinitesimal length element $ds$ at an arbitrary point of the first surface and the corresponding element $ds'$ at the image point does not depend on the chosen direction at the first point. He recalls that this condition is the one that Lambert, Lagrange and Gauss adopted as a general principle in their theory of geographic maps. In  Notes V and VI to Monge's treatise, Liouville gives his own solution to this problem.

Let us now say a few words on the work of Chebyshev on cartography.

 Chebyshev (1821--1894) was a devoted reader of Lagrange.\footnote{Chebyshev's  biographer in \cite{bio} reports that Chebyshev thoroughly studied the works of Euler, Lagrange, Gauss, Abel, and other great mathematicians. He also writes that, in general, Chebyshev was not interested in reading the mathematical works of his contemporaries, explaining that spending time on that would prevent him of having original ideas.} Like Euler, he was interested in almost all branches of pure and applied mathematics, and he spent a substantial part of his time working on industrial machines  of all kinds. He conceived a computing machine, a walking machine, and other kinds of machines, and he worked on linkages and hinge mechanisms from the practical and theoretical points of view; cf. \cite{Sossinsky}. Approximation theory and optimization were among his favorite fields. Thus, it is not surprising that the questions raised by the drawing of geographic maps naturally attracted his attention.  

Chebyshev's \emph{Collected papers} \cite{T-oeuvres} contain two papers on cartography, \cite{Cheb1} and \cite{Cheb2}, both called \emph{Sur la construction des cartes g\'eographiques}. The two papers contain several interesting ideas on this subject. 

In the paper \cite{Cheb1}, Chebyshev starts by recalling that it is easy to reproduce an arbitrary part of the globe while preserving angles (he writes: ``such that there is constantly a similarity between the infinitely small elements and their representation on the map.") The problem, he says, is that \emph{the magnification ratio} (a kind of scaling factor) in such a map, varies from point to point. This is one of the reasons for which the representation of the curved surface that we see on the map is not faithful -- it is deformed. The \emph{scale} of the map is not constant, it depends on the chosen point. In the general case, the magnification ratio depends on the choice of a point together with a direction at a point. For mappings that are angle-preserving (which Chebyshev call  ``similarities at the infinitesimal level"), the magnification ratio depends only on the point on the sphere and not on the chosen direction. Thus, the magnification ratio is a function defined on the sphere.

Chebyshev studies then the problem of finding geographical maps that preserve angles and such that this magnification ratio is minimal. To solve this problem, he starts by a formula of Lagrange from his paper \cite{Lagrange1779} for the magnification ratio of such maps (Formula  (\ref{distortion}) again), which he writes in the following form:
  \[m=\frac{\sqrt{f'(u+t\sqrt{-1})F'(u-t\sqrt{-1})}}{\displaystyle\frac{2}{e^u+e^{-u}}},
\] 
using the notation of Lagrange which we already recalled, and where $f$ and $F$ are again arbitrary functions.
The formula leads to 
\[\log m = \frac{1}{2}\log \big( f'(u+t\sqrt{-1})\big) +  \frac{1}{2}\log \big( F'(u-t\sqrt{-1})\big) -\log \frac{2}{e^u+e^{-u}}.\]
The first two terms in the right hand side of the last equation, which contain arbitrary functions, is the solution of the Laplace equation:
\[\frac{\partial^2 U}{du^2}+\frac{\partial^2 U}{dt^2}=0.
\]

From the theory of the Laplace equation, Chebyshev concludes that the \emph{minimum} of the deviation of a solution of this equation from the function $\displaystyle \log \frac{2}{e^u+e^{-u}}$, on a region bounded by an arbitrary simple closed curve is attained if the difference 
\[U-\log \frac{2}{e^u+e^{-u}}\]
is constant on the curve.
Thus, one obtains the value of 
\[U= \frac{1}{2}\log \big( f'(u+t\sqrt{-1})\big) +  \frac{1}{2}\log \big( F'(u-t\sqrt{-1})\big)\]
up to a constant, and from there, the value of the functions
\[ f'(u+t\sqrt{-1})\]
and 
\[ F'(u-t\sqrt{-1}).\]
Up to a constant factor, the functions that give the best projection are found.

 Chebyshev then considers the important special case studied by Lagrange, where the images of the parallels and meridians by the projection map are circles in the plane.  In this case, Lagrange's formula for the magnification ratio takes the form:
 \[m=\frac{1}{\frac{2}{e^u+e^{-u}}\big(a^2e^{2cu}+2ab\cos 2c(t-g)+b^2e^{-2cu}
 \big)
 }.
\]
Chebyshev then works out the details of the result saying that the best geographical mapping is obtained when $m$ is constant on the curve. He obtains precise estimates in the case where the curve bounding the region which is represented is close to an ellipse. He concludes the first paper by highlighting the following general result:
\begin{quote} \emph{There exists an intimate connection between the form of a country and its best projection exponent}.
\end{quote}
More details on this question are given in his second paper \cite{Cheb2}.

 The second paper  is a sequel to the first. In the introduction, Chebyshev, who was particularly interested in practical application of mathematics, places  this problem in a most general setting, namely, that of ``finding the way to the most advantageous solution of a given problem." Going into more detail regarding this question would take us too far. Let us only recall that the theory of optimization was one of the favorite subjects of Chebyshev.
 
 Returning to the problem of geographical maps, Chebyshev places it in the setting of the calculus of variations. 
 He uses the computations made in his first paper \cite{Cheb1}, and in particular the relation with the Laplace equation.
  He returns to the question of the dependence of the choice of the most appropriate mapping on the shape of a country, that is, on the form of its boundary and to his results deduced from the work of Lagrange saying that the best projection is one for which the ratio of magnification is constant on the boundary of the country.

   The details of Chebyshev's arguments are somehow difficult to follow, even if the general ideas are clear. Milnor, in his paper  \cite{Milnor}, gives an exposition of the same result with complete proofs following Chebyshev's ideas. 
  We report briefly on Milnor's exposition.
  
  Let $\Omega$ be a simply connected open subset of the sphere of radius $r$ in Euclidean 3-space. We consider a conformal mapping $f$ from $\Omega$ onto the Euclidean plane $E$. Such a map has a well-defined \emph{infinitesimal scale} at each point $x$ in $\Omega$, defined as 
  \[\sigma(x)=\lim_{y\to x} \frac{d_E(f(x),f(y))}{d_S(x,y)}.\]
  The limit exists because the mapping $f$ is conformal (the deviation at a point does not depend on the direction).
  
 For such a conformal map $f$, its  \emph{distortion} is defined as the ratio
 \[\frac{\sup_{x\in \Omega}\sigma(x)}{\inf_{x\in \Omega}\sigma(x)}.\]
 
The goal is to find a map with the smallest distortion. We have the following result:
 \begin{quote}
 [Chebyshev's theorem -- Milnor's revision] \emph{Assume that the boundary of the open set $\Omega$ is a twice differentiable curve. Then there exists one and (up to similarity) only one conformal mapping $f$ from $\Omega$ onto the Euclidean plane with minimal distortion. This map is characterized by the property that its infinitesimal scale function $\sigma(x)$ is constant at the boundary.}
 \end{quote} 
   Milnor adds:
\begin{quote}\small
This result has been available for more than a hundred years, but to my knowledge it has never been used by actual map makers.
\end{quote}
 Milnor notes after his proof of this result that the ``best possible" conformal map $f$ that is given by this theorem is locally injective, but in general it is not globally injective. He adds that $f$ will be globally injective if $U$ is geodesically convex, and that in that case $f(U)$ is also geodesically convex. 

In practice, very few countries are geodesically convex, but the result may be used to draw maps of geodesically convex regions of the Earth.

 Let us note by the way that the definition of the infinitesimal scale function, and especially its global version which is also considered by Milnor, is very close to some basic definitions that appear in Thurston's theory of stretch maps between surfaces \cite{Thurston1985}.  The question of finding the best geographical map between surfaces bears a lot of analogy to that of finding the best Lipschitz map in the sense of Thurston.

In Milnor's exposition of Chebyshev's ideas, the search for the best geographical map amounts to the solution of the Laplace equation with given boundary values. Potential theory thus appears at the forefront. 

In the same paper, Milnor continues his cartographical investigations by  studying the magnification ratio for various kinds of geographical maps. He obtains a lower bound for this ratio in the case of the so-called azimuthal equidistant projection, and he proves the existence and uniqueness of a minimum distortion map in this case.

Darboux in  \cite{Darboux-Chebyshev} also gave a solution to the same problem, which is based on Chebyshev's ideas and using potential theory. It seems that Milnor was not aware of the work of Darboux on this question, since he does not mention it.

Summarizing the contribution of Chebyshev, we can highlight the following three reasons for which that work, and the works of the other cartographers that we mentioned, are interesting from the point of view of the subject of our present survey:
  
\begin{enumerate}

\item The work involves definitions and notions related to mappings between surfaces which are close to those considered by quasiconformal theorists;

\item The techniques used in cartography are those of partial differential equations, as inaugurated by Euler and used later on by quasiconformal mappers;

\item  Chebyshev proved existence and uniqueness results on extremal mappings under some conditions.

\end{enumerate}

There is more, namely, ideas related to \emph{moduli}. Chebyshev, in his paper \cite{Cheb2} (p. 243), recalls that there are several kinds of angle-preserving maps from the sphere to the plane. On p. 244ff. of that paper, he addresses the question of how the map used changes, \emph{when the boundary of the region to be represented varies}. He discusses the variation of a point called the \emph{center} and the exponent of the projection, and he mentions in particular the cases of regions bounded by pieces of second-degree curves.  

For a more detialed exposition on Chebyshev's work on geography and relations with his other works, the reader is referred to the paper \cite{2016-Tchebyshev}.

 One of the last prominent representatives of the differential-geometric work on geography is Beltrami. He worked on problems related to cartography in the tradition of Gauss. One of his important papers in this field is the \emph{Risoluzione del problema: Riportare i punti di una superficie sopra un piano in modo che le linee geodetiche vengano rappresentate da linee rette} (Solution of the problem: to send the points of a surface onto a plane in such a way that the geodesic lines are represented by straight lines)\footnote{This is the paper which contains Beltrami's famous result saying that a Riemannian metric on a surface which can be locally mapped onto the plane in such a way that the   geodesics are sent to Euclidean lines has necessarily constant curvature.} (1865)   \cite{Beltrami1865}. Beltrami declares in this paper that a large part of the research done before him on similar questions was concerned with conservation either of angles or of area. He says that even though these two properties are considered as the simplest and most important ones for geographical maps, there are other properties that one might want to preserve. He declares that since the projection maps that are used in geography are  mainly concerned with measurements of distances, one would like to exclude maps where the images of distance-minimizing curves are too remote from straight lines. Thus, we encounter again a  notion of a map sending a geodesic to a quasigeodesic. (We already mentioned that this property was considered in a paper by Euler, \cite{Euler-pro-Desli-1777}.) Beltrami calls the maps between surfaces that he considers ``transfer maps."  He mentions incidentally that the central projection of the sphere is the only map that transforms the geodesics of the sphere into Euclidean straight lines.  He declares that beyond its applications to the drawing of geographical maps, the problem may lead to ``a new method of geodesic calculus, in which the questions concerning geodesic triangles on surfaces can all be reduced to simple questions of plane trigonometry."  This is a quite explicit statement in which one can see that questions on geography acted as a motivation for the development  of the field of differential geometry of surfaces. There are other such explicit statements. For instance, Darboux gave a talk at the Rome 1908 ICM talk whose title is \emph{Les origines, les m\'ethodes et les probl\`emes de la g\'eom\'etrie infinit\'esimale} (The origins, methods and problems of infinitesimal geometry) \cite{Darboux-ICM}. He declares there that the origin of infinitesimal geometry lies in cartography. He starts by an exposition of the history of the subject. We quote an interesting passage:\footnote{In the present article, the translations rom the French are mine.}
\begin{quote} \small 
Like many other branches of human knowledge, infinitesimal geometry was born in the study of practical problems. The Ancients were already busy in obtaining plane representations of the various parts of the Earth, and they had adopted the idea, which was so natural, of projecting onto a plane the surface of our globe. During a very long period of time, people were exclusively attached to these methods of projection, restricting simply to the study of the best ways to choose, in each case, the point of view and the plane of the projection. It was one of the most penetrating geometers, Lambert, the very estimated colleague of Lagrange at the Berlin Academy, who, pointing out for the first time a property which is common to the Mercator maps and called \emph{reduced maps} and to those which are provided by the stereographic projection, was the first to conceive the theory of geographical maps from a really general point of view. He proposed, with all its scope, the problem of representing the surface of the Earth on a plane with the similarity of the infinitely small elements. This beautiful question, which gave rise to the researches of Lambert himself, of Euler, and to two very important memoirs of Lagrange, was treated for the first time in all generality by Gauss. [...] Among the essential notions introduced by Gauss, one has to note the systematic use of the curvilinear coordinates on a surface, the idea of considering a surface like a flexible and inextensible fabric, which led the great geometer to his celebrated theorem on the invariance of total curvature, to the beautiful properties of geodesic lines and their orthogonal trajectories, to the generalization of the theorem of Albert Girard on the area of the spherical triangle, to all these concrete and final truths which, like many other results due to the genius of the great geometer, were meant to preserve, across the ages, the name and the memory of the one who was the first to discover them.
\end{quote} 
\section{The work of Tissot}\label{s:Tissot}
 
In this section, we review the work of Tissot. This mathematician-geographer discovered new projections of the sphere that are useful in cartography, some of which were still used recently. But most of all, Tissot brought substantial new theoretical ideas in this field. In particular, in his paper \cite{Tissot1881}, he developed the analytical tools needed to obtain formulae for the minimum of the distortion. His most important contribution, which is certainly the reason for which Teichm\"uller quotes him in his paper \cite{T20}, is the discovery of a device that survived under the name \emph{Tissot indicatrix}, or \emph{ellipse indicatrix}. As we shall see, the underlying idea is at the basis of the notion of quasiconformal mapping which were later developed by Gr\"otzsch, Lavrentieff, and Teichm\"uller. The Tissot indicatrix is also mentioned explicitly and several times in the paper \cite{Groetzsch1930} by Gr\"otzsch. We shall discuss this paper in the section on G\"otzsch, in \S\ref{s:G} of the present paper.

The Tissot indicatrix was used during more than a hundred years in the drawing of geographical maps. Let us summarize right away the idea. 

We consider a mapping between two surfaces. For each point of the domain surface, we take an infinitesimal circle centered at that point. If one does not want to talk about infinitesimals, he may consider that such a circle is in the tangent space at the given point, centered at the origin. In this way, the domain surface is equipped with a field of infinitesimal circles. If the mapping is conformal, then it sends infinitesimal circles to infinitesimal circles. But the ratios of the radii (Tissot says: ``to the extent that one can talk about radii of infinitesimal circles") of the image circles to those of the original ones vary from point to point. For what concerns us here, we may assume that the radii of the infinitesimal circles of the first surface are  constant. The collection of ratios obtained in this manner constitutes a measure of the amount by which the mapping distorts distances. In the general case, the image of an infinitesimal circle is an infinitesimal ellipse which is not necessarily a circle. The collection of shapes of the image ellipses (that is,  the ratio of their major to the minor axes) is another measure of the distortion of the mapping.  An area-preserving mapping is one for which the image ellipses have the same area as the original circles. The inclination of the ellipse (which is defined, like the ratio of the great axis to the small axis, provided the ellipse is not a circle), measured by the angle its great axis makes with the horizontal axis of the plane, is another measure of the distortion of the mapping. Knowing the axes of the  ellipse, the angle distortion of the mapping is known. Thus, the study of the field of ellipses gives information on: 
\begin{enumerate}
\item the non-conformality of the map, that is, its angle distortion;
\item the area distortion;
\item the distance distortion.
\end{enumerate}
By setting the radii of the infinitesimal circles on the domain surface to be equal to one, then the shape,  inclination and area of each image ellipse, which Tissot says is ``a kind of \emph{indicatrix}," is expressed by some finite data which give precise information on the distorsion of the map. 

In his work,   Tissot establishes analytical formulae for the lengths and the directions of the axes of the image ellipses.
He also considers the \emph{maximum} of the distortion over the surface. Similarly, he studies the ratios according to which the length element is changed in the various directions, the greatest and the smallest such ratios (which are precisely equal to the semi-axes of the ellipses), and finally, the change in area. 

\begin{figure}[htbp]

\centering
\includegraphics[width=12cm]{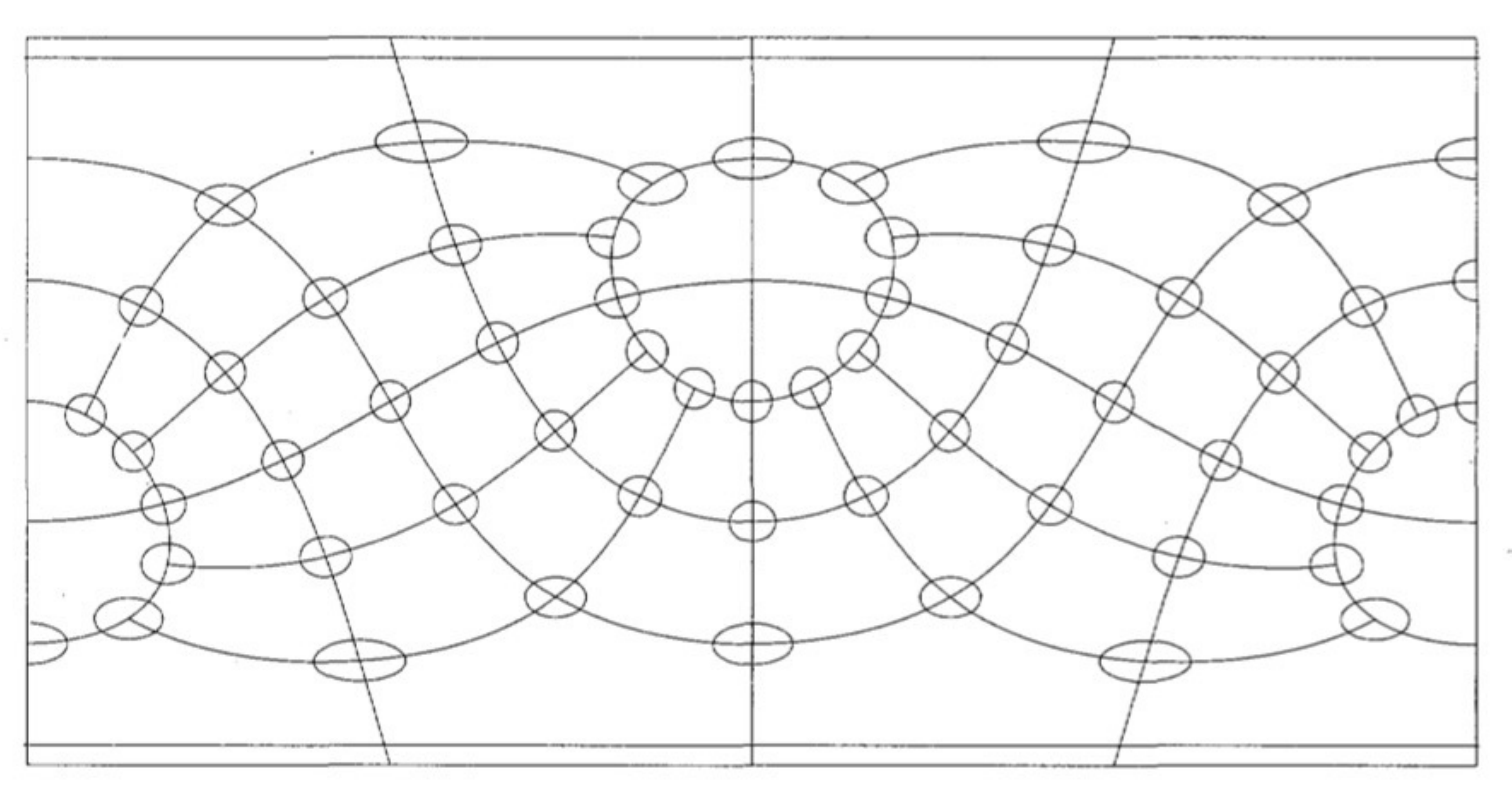}
\caption{\smaller{The Tissot indicatrix for the Oblique Plate Carr\'ee projection; from the \emph{Album of map projections} \cite{Snyder-Album} p. 28.}} \label{Tissot-Oblique}
\end{figure}

Cartographers used to draw the \emph{Tissot indicatrix}, that is, the field of ellipses of a geographical map. In Figure \ref{Tissot-Oblique} is represented the field of ellipses corresponding to the so-called \emph{oblique plate carr\'ee} (also called equirectangular) projection. In such a projection, the meridians and parallels form a pair of perpendicular foliations. Tissot considers this foliation as a grid dividing the surface into equal infinitesimal rectangles whose side ratios are equal to some fixed ratio between meridian and parallel distances.

Let us emphasize the fact that the Tissot indicatrix is not a measure of conformality in the strict sense, since it involves distances and areas and not only the conformal structures of the surfaces. It is however a measure of conformality when the given conformal structure underlies a metric structure. We shall talk more about this below, when we consider quasiconformal mappings, in particular when we talk about the so-called length-area method.
 
Tissot first published his work in a series of Comptes Rendus notes (see \cite{Tissot-CR1}--\cite{Tissot-CR4}.) Several years later, he gave the details of his work in installments (Tissot \cite{Tissot1878}--\cite{Tissot1880}), and then in the long memoir \cite{Tissot1881}, which reproduces the articles \cite{Tissot1878}--\cite{Tissot1880} and additional material. 
The paper \cite{Tissot1878a} and the introduction of \cite{Tissot1881} contain an outline of the whole theory.  

Tissot worked in the setting of  arbitrary differentiable maps between surfaces, although primarily in the case where the domain surface is a subset of a sphere or spheroid, and the range a subset of the Euclidean plane. At several places he considers ``any map (which he calls a \emph{representation}) from a surface onto another." On p. 52 of his paper \cite{Tissot1878}, he mentions applications of his theory to the case where the domain is an arbitrary  surface of revolution. We recall that this special case   was also thoroughly investigated by Euler. Tissot proves the following fundamental (although easy to prove) lemma (\cite{Tissot1878} p. 49 and \cite{Tissot1881} p. 13):
  
  \begin{quote}
   \emph{At every point of the surface which has to be represented, there are two mutually perpendicular tangents, and if the angles are not preserved then there is a unique such pair, such that the corresponding tangent vectors on the other surface also intersect at right angles.} 
   \end{quote}

   Let us insist on the fact that with the above condition of not preserving angles, the setting of Tissot's study becomes different from that of his predecessors (in particular Lagrange and Chebyshev) who studied mainly angle-preserving maps with some other distortion parameters minimized.
   
   From this lemma, Tissot deduces the following:
  
    \begin{quote}
     \emph{On each of the two surfaces there exists a ``system of orthogonal trajectories," and if the representation is nowhere angle-preserving, there exist only such one pair, which are preserved by the map.}
     \end{quote}
        
In modern language, each system of trajectories in the pair is a foliation of the surface. It is interesting to see how Tissot expresses his ideas without the specialized geometric terminology (which was still inexistent). For the network consisting of two perpendicular foliations, He uses the word ``canevas" (with reference to a canvas fabric). 

From the practical point of view, these orthogonal foliations provide a system of ``Euclidean" coordinates on the surface.\footnote{The existence of two orthogonal foliations is also an important feature of the famous Mercator map, conceived by the mathematician and geographer Gerardus Mercator (1512--1594). Mercator was born under the name Geert (or Gerhard) De Kremer, meaning shopkeeper in Flemish. (Mercator is the Latinized version.)  The Mercator map is conformal and it is based on a cylindrical projection of the sphere in which meridians and parallels become orthogonal straight lines. Mercator's maps were  extensively used by navigators since the sixteenth century.  Similar maps were known in antiquity. They were used since the work of Eratosthenes (third century B.C.), in particular to represent the new lands conquered by Alexander the Great.} Tissot says that the surface admits infinitely many such networks, but that those associated to a map (projection) between two surfaces are unique. He formulates the following rule satisfied by an arbitrary deformation, which is independent of the nature of the surface or of the given map: 
   
   \begin{quote}
   \emph{Any map between two surfaces may be replaced, around each point, by an orthogonal projection performed at a convenient scale.}
   \end{quote}
   
    The rule is proved in a geometric way using the lemma we mentioned above. Tissot deduces from this rule a large number of properties.

 In the second chapter of the memoir \cite{Tissot1881}, Tissot  reports on specific applications of his theory.
 One important idea in this chapter is formulated right at the beginning of the memoir: 
 \begin{quote} 
 \emph{To find the projection mode that is the most appropriate to the plane representation of a given country.}
 \end{quote}  
 In other words, the projection mode will depend on the shape of the country. Thus, the emphasis is, like in Chebyshev's work, on the boundary of the region of the sphere (or the spheroid) that has to be represented. 

 In the memoir \cite{Tissot1881}, Tissot studies thoroughly projections that produce very small angle distortion (of the order of a few seconds), and where the length distortions are minimized. He says that a similar method may be applied to projections which distort areas by negligible quantities, and where angle distortion is minimized.

Tissot's memoir \cite{Tissot1881} also contains a section on the history of cartography (p. 139ff.).

  We refer the reader to the article \cite{2016-Tissot} for a more extended survey of the work of Tissot.
  
Darboux, who, as we already mentioned, was also interested in geography, wrote a paper on Tissot's work \cite{Darboux-Tissot}.  In the introduction, he recalls  the techniques introduced by Tissot for the representation of a given country, and in particular his use of power series expansions of functions, but he says that ``[Tissot's] exposition appeared to me a little bit confused, and it seems to me that while we can stay in the same vein, we can follow the following method ..." He then explains in his own way Tissot's work, making the relation with the works of Gauss, Chebyshev and Beltrami on the drawing of geographical maps. In particular, he extends Tissot's theory to maps between surfaces considered in the differential-geometric setting of Gauss's theory, starting with  arbitrary curvilinear coordinates with length element
   \[ds^2=Edu^2+2Fdudv+Gdv^2\]
   and using the theory of conformal representations.

We already mentioned another paper by Darboux on geographical maps, \cite{Darboux-Chebyshev}, in which he  gives a detailed proof of the result of Chebyshev saying that the most advantageous representation of a region of the sphere onto the Euclidean plane is the one where the magnification ratio is constant on the boundary of the surface to be represented. 

Let us also mention a third paper that Darboux wrote on cartography, \cite{Darboux-Lagrange}, in which he solves a problem addressed by Lagrange in the paper \cite{Lagrange1779} that we already mentioned. The question concerns a quantity which Lagrange calls the ``exponent of the projection." This problem is reduced to the following question in elementary geometry:
\begin{quote}
 \emph{Given three points on the sphere, can we draw a geographical map, with a given exponent, such that these three points are represented by three arbitrarily chosen points on the map?} 
 \end{quote}
 Lagrange declares in his paper that a geometric solution to this problem seems very difficult, and that he did not try to find a solution using algebra. Darboux solves the problem in a geometric manner. He declares that it is the recent progress in geometry that made this solution possible.

    One can say much more about the contributions of geographers on close-to-conformal mappings, but we must talk now about the modern period.

In the next three sections, we briefly report on the lives and works of four pioneers of the modern theory of quasiconformal mappings: Gr\"otzsch, Lavrentieff, Ahlfors and Teichm\"uller.   The questions addressed become, in comparison with those that were asked by geographers, as follows:
\begin{enumerate}
\item \label{q1} The desired map is a close-to-conformal mapping between two Riemann surfaces (that is, surfaces equipped with a conformal structures). The metric structure is not given in the formulation of the problem, although it is used as a tool in the developments.  

\item The two surfaces are no more restricted to be subsets of the sphere and of the Euclidean plane respectively. We already saw that several authors already worked in this generality (Euler, Lagrange, Chebyshev, etc.), although the surfaces they considered were simply connected, or at least planar (homeomorphic to subsets of the plane).

\item \label{q2} The desired extremal map minimizes the supremum over the surface of the ratios of the big axes to the small axes of the infinitesimal Tissot ellipses.

\item \label{q3} There are topological restrictions on the map. For instance, in the case of mappings between closed surfaces,  one asks that they are in a given homotopy class. The theory also involves the case of mappings that preserve distinguished points. For surfaces with boundary, these distinguished points may be on the boundary. In the theory of quasiconformal mappings of the disc, the homotopies defining the equivalence relation are required to be the identity on the boundary of the disc. There are several other variants of the theory.
\end{enumerate}

Let us make a comment on (\ref{q1}). In principle, the theory of quasiconformal mappings does not involve distances, since, being maps between Riemann surfaces, they are sensible only to the conformal structure, that is, the angle structure on the surface, and how it is distorted. However, when the conformal structure is induced by a metric, relations between the angle structure and length and area arise, and we are led to questions similar to those in which cartographers were interested in length and area distortion of geographical maps. The ``length-area method," was one of the essential tools in the works of Gr\"otzsch, Ahlfors and Teichm\"uller on quasiconformal mappings. We shall return to this question in the second part of our paper. We also mention Gr\"otzsch's paper \cite{Groetzsch1930} which contains precise estimates on length distortion of quasiconformal mappings between multiply-connected domains of the plane.

The consideration of distinguished points mentioned in (\ref{q3}) is essential. Typically, if the surface is a disc, then the Riemann Mapping Theorem says that (in the case where both surfaces are not conformally equivalent to the plane) one can find a \emph{conformal} mapping between them, and that furthermore there exists such a mapping with the property that any chosen three pair of distinct distinguished points on the boundary of the first surface are sent to any chosen three distinct points on the boundary of the second surface. The next step starts with discs with \emph{four} distinguished points on their boundary, and this is the case considered by Gr\"otzsch. 

\part{Quasiconformal mappings}
We start the second part of our paper with the work of Gr\"otzsch. Then we continue with the works of Lavrentieff, Ahlfors and Teichm\"uller.
  \section{Gr\"otzsch} \label{s:G}

  Camillo Herbert Gr\"otzsch (1902--1993) studied mathematics at the University of Jena, between 1922 and 1926, where Paul Koebe was one of his teachers. Koebe moved to Leipzig  in 1926, and Gr\"otszch followed him. He worked there on his doctoral dissertation, under the guidance of Koebe, and obtained his doctorate in 1929, cf. \cite{Groetzsch-Diss}.
  Reiner K\"uhnau, who was Gr\"otzsch's student, wrote two long surveys on his work, \cite{Kuhnau1} and \cite{Kuhnau2}. Between the years 1928 and 1932, Gr\"otzsch published 17 papers on conformal and quasiconformal mappings. He introduced the notion of extremal quasiconformal mapping, which he called ``m\"oglichst konform" (conformal as much as possible). Gr\"otzsch's papers are generally short (most of them have less than 10 pages) and they almost all appeared in the same journal, the \emph{Leipziger Berichte}.\footnote{This journal is sometimes referred to under other names, see Gr\"otzsch's references in the bibliography of the present paper.} In fact, Gr\"otzsch published almost all his papers in rather unknown journals. 
   K\"uhnau writes in \cite{Kuehnau-H} that Gr\"otzsch ``did not like to be dependent on the grace of a referee. For him only Koebe was a real authority." His papers are quoted in the literature on quasiconformal mappings, but it seems that nobody read them. Volume VII of the \emph{Handbook of Teichm\"uller theory} will contain translations of several of these papers (cf. \cite{Gr1928}, \cite{Groetzsch1930}, and there are others). K\"uhnau adds in \cite{Kuehnau-H} that ``to his disadvantage was also the fact that Gr\"otzsch almost never appeared at a conference. Ahlfors told me that he earlier thought that Gr\"otzsch did not really exist." 
 
  Gr\"otzsch is mostly known for his solution of a problem which is easy to state, namely, to find the best quasiconformal mapping between two quadrilaterals  (discs with four distinguished points on their boundary). According to Ahlfors (\cite{Ahlfors1964} p. 153), this problem was first considered as a mere curiosity. Things changed, and this work is now considered as a key result, especially since the generalization of this problem and it solution by Teichm\"uller. Talking about the method of Gr\"otzsch, Ahlfors writes in his 1953 paper celebrating the hundredth anniversary of Riemann's inaugural dissertation  \cite{Ahlfors-dev}:
\begin{quote}\small
Comparison between length and area in conformal mapping, and the obvious connection derived from the Schwarz inequality, has been used before, notably by Hurwitz and Courant. The first to make systematic use of this relation was H. Gr\"otzsch, a pupil of Koebe. The speaker hit upon the same method independently of Gr\"otzsch and may, unwillingly, have detracted some of the credit that is his due. Actually, Gr\"otzsch had a more sophisticated point of view, but one which did not immediately pay off in the form of simple results.
\end{quote}
  
  In fact,  Gr\"otzsch studied, besides quasiconformal mappings between rectangles, maps between surfaces that are not necessarily simply-connected. For instance, his paper \cite{Groetzsch1930} concerns quasiconformal mappings between multiply-connected subsets of the plane. The problem of finding the \emph{best} quasiconformal mapping between two homeomorphic Riemann surfaces is not addressed (this was one of Teichm\"uller's achievements). Rather, in this paper, Gr\"otzsch obtains estimates for the distortion of length and area under a quasiconformal mapping.  It is interesting to see that Gr\"otzsch, in the paper \cite{Groetzsch1930}, does not use any word equivalent to the word ``quasiconformal" to denote these maps. In fact, he does not give them a name but only a notation, and he extensively uses the setting of Tissot's work. Let us be more precise.
  
 Gr\"otzsch considers in \cite{Groetzsch1930} bijective and continuous mappings of a domain $B$ that are uniform limits of affine mappings except for a finite number of points. This implies a differentiablity property (which Gr\"otzsch does not state explicitly). He then assumes that at each point where the map is differentiable, the ratio $a/b$ of the great axis to the small axis of an infinitesimal ellipse which is the image of an infinitesimal circle satisfies the condition 
  \[1/Q\leq a/b\leq Q\]
  for some $Q\geq 1$ that is independent of the choice of the point. Gr\"otszch formulates this result as follows:
  \begin{quote}
  \emph{The distorsion of the shape of Tissot's indicatrix (distortion ellipse) of the mapping is generally kept within fixed bounds.}
  \end{quote}
  
  Gr\"otzsch then says that a mapping $B$ satisfying these conditions is called a mapping $\frak{A}_Q$ of $B$.
   This suggests that the map is $Q$-quasiconformal (in the modern sense), but no special name is used for such a map.

  Teichm\"uller mentions extensively the papers of Gr\"otzsch, especially in his later papers.  In the introduction to his paper \cite{T23}, he writes ``Much of what is discussed here is already contained in the works of Gr\"otzsch, but mostly hidden or specialized in typical cases and in a different terminology."
  
   Apart from his method for the modulus of a rectangle which is repeated in several subsequent papers and books, Gr\"otzsch's results on quasiconformal mappings are still very poorly known.  Later on, Gr\"otzsch became known for his work on graph theory.  The reader is referred to the article \cite{Kuehnau-H} by K\"uhnau for a lively information on Gr\"otzsch in relation with Teichm\"uller. 
  
\section{Lavrentieff} \label{s:L} 
Unlike Gr\"otzsch, Mikha\"\i l Alekse\"\i evitch Lavrentieff\footnote{Like all the Russian similar names, Lavrentieff may also be written with a terminal v instead of ff. We are following Lavrentieff's own transcription of his name in his papers written in French, cf. \cite{Lavrentieff1927a} \cite{Lavrentieff1927b} \cite{Lavrentiev-ICM}  \cite{Lavrentiev1935a} \cite{Lavrentiev1935} \cite{Lavrentieff-Actualites}.}  (1900--1980) had a brilliant career, as a mathematician and as a leader in science organisation. 
    
    Between 1922 and 1926, Lavrentieff worked on his doctoral thesis, on topology and set theory,\footnote{At the epoch of Luzin, the expressions ``topology" and ``set theory" often indicated the same subject.  The reader may recall that one of the first books on modern topology is Hausdorff's \emph{Grundz\"uge der Mengenlehre} (Foundations of set theory) (1914). The second edition was simply called \emph{Mengenlehre}.} with Luzin as an advisor.  He later remained faithful to the Luzin school despite the terrible campaign that was made around him.\footnote{Nikolai Luzin (1883-1950) is the main founder of the Moscow school of  function theory, and one of the architects of descriptive set theory.  In 1936, after a very successful career in mathematics, Luzin fell into disfavor, and the Soviet authorities started a violent political campaign against him. Six years before, his  mentor, Egorov, had experienced a similar disgrace, due to his religious sympathies. Egorov  died miserably in detention, from a hunger strike, in 1931, after having been dismissed in 1929 from his position at the university, and jailed in 1930 with the accusation of being a ``religious sectarian." Luzin lost his position at the university. His case is famous, and what became known as the ``Luzin affair" shook the scientific world in the Soviet Union and beyond. The story is recounted in several books and articles, e.g. \cite{Demidov}.   In 1974,  Lavrentieff published a memorial article on Luzin \cite{L1974}.  The book \cite{Lavrentieff-century} which appeared on the occasion of Lavrentieff's hundredth anniversary \cite{Lavrentieff-century} contains a chapter on the Luzin affair. The relation between Luzin and the Lavrentieff family remained constant.}

      In 1927, Mikha\"\i l  Lavrentieff made a stay in Paris during which he followed courses and seminars by Montel, Borel, Julia, Lebesgue and Hadamard. At this stage his interests shifted from topology to the theory of functions of a complex variable.  He published two Comptes Rendus notes on the subject, \cite{Lavrentieff1927a} and \cite{Lavrentieff1927b}. The first one concerns boundary correspondence of conformal representations of simply connected regions. In this paper, Lavrentieff studies the \emph{ratios of the lengths} of the boundaries of the two discs. We recall that the same kind of problem was addressed by Chebyshev in his study of geographical maps. Among the results that Lavrentieff obtained in the paper \cite{Lavrentieff1927a} are the following two:
      
      \begin{quote}\small
      \emph{1.--- Let $D$ and $D'$ be two domains in the plane, bounded by two simple closed curves $\Gamma$ and $\Gamma'$ whose curvature is bounded. Consider a conformal mapping between $D$ and $D'$. Then the ratio of lengths of corresponding arcs of $\Gamma$ and $\Gamma'$ is bounded.
      }
      \end{quote}
         \begin{quote} \small
         2.--- Let $D$ and $D'$ be two domains in the plane, bounded by two simple closed curves $\Gamma$ and $\Gamma'$ which have continuously varying tangents. Consider a conformal mapping between $D$ and $D'$. If $\delta$ and $\delta'$ denote the the lengths of corresponding arcs on $\Gamma$ and $\Gamma'$, then we have 
         \[K_1\delta^{1-\epsilon}>\delta' > K_2\delta^{1+\epsilon}
         \]
         where $K_1$ and $K_2$ are constants that depend only on $\epsilon$.
          \end{quote}
         
      The second paper, \cite{Lavrentieff1927b}, concerns function theory.

Lavrentieff founded new areas of research in mechanics and he applied  the theory of functions of a complex variable to the study of non-linear waves. His interest in quasiconformal mappings was not independent of his work on physics.

  Some of Lavrentieff's papers are written in French, and his works were known in the West.  Courant, in his 1950 book \cite{Courant1950},  has a section on a method developed by Lavrentieff to show that some extremal domains have analytic boundaries.  Lavrentieff did important work in  aerohydrodynamics which led him very far in the applications. In particular, he conceived a new kind of bullet, known under the name Katyusha, which played a very important role during World War II. It is possible that his work on quasiconformal mappings was of great help in that domain. 
Lavrentieff received the most prestigious prizes awarded in the Soviet Union, among them the Lenin Prize and the Lomonosov gold medal.

  Lavrentieff's work on quasiconformal mappings dates back to his early career. In 1928, he gave a talk at the Bologna ICM \cite{Lavrentiev-ICM} in which he described an application of mappings which in some specific sense are quasiconformal mappings. The question was to construct the Riemann Mapping Theorem using a sequence of explicit mappings obtained from the theory of partial differential equations, using a minimization principle. We note incidentally that a similar application of quasiconformal mappings is mentioned by Teichm\"uller in the last part of his paper \cite{T20}. In his Comptes Rendus note \cite{Lavrentiev1935a}  and his paper \cite{Lavrentiev1935}, both published in 1935, Lavrentieff returns to this subject with more details. He introduces the notion of quasiconfomal mappings with a very  weak condition of differentiablilty. The two papers are titled ``Sur une classe de repr\'esentations \emph{continues}." Let us give a few more details.
  
  The Note \cite{Lavrentiev1935a} starts with the definition of an \emph{almost analytic} mapping (``fonction presque analytique") $f$.
This definition involves two functions, $p$ and $\theta$, which he calls the \emph{characteristic functions} of $f$. The function $p$ is the ratio of the great axis to the small axis of the ellipse which is the image by $f$ of an infinitesimal circle. The function $\theta$ is the angle between the great axis of that ellipse and the real axis.  The name  \emph{characteristic function} reminds us of the Tissot characteristic which we mentioned in \S \ref{s:Tissot}. The definitions introduced by Tissot and Lavrentieff are very close. Lavrentieff was probably not aware of the work of Tissot. The definition that Lavrentieff gives in \cite{Lavrentiev1935a}  and (\cite{Lavrentieff1927b} p. 407) is the following:

A function $w=f(z)$ of a complex variable $z$ in a domain $D$ is termed \emph{almost analytic} if the following three properties hold:
\begin{enumerate}
\item $f$ is continuous on $D$.

\item In the complement of a countable and closed subset of $D$, $f$ is an orientation-preserving local homeomorphism.

\item There exist  two functions, $p(z)\geq 1$ and $\theta(z)$ defined on $D$ such that
\begin{itemize}
\item In the complement of a subset $E$ of $D$ consisting of a finite number of analytic arcs, $p$ is continuous and $\theta$ is continuous at any point $z$ satisfying $p(z)\not=1$.

\item On any domain $\Delta$ which does not contain any point of $E$ and whose frontier is a simple analytic curve, $p$ is uniformly continuous, and if such a domain $\Delta$ and its frontier do not contain any point such that $p(z)=1$, then $\theta$ is also uniformly continuous on $\Delta$.

\item For any point $z_0$ in $D$ which is not in $E$, let $\mathcal{E}$ be the ellipse centered at $z_0$ such that the angle between the great axis of $\mathcal{E}$ and the real axis of the complex plane is $\theta(z)$ and $p(z)\geq 1$ the ratio of the great axis to the small axis of $\mathcal{E}$.
Let $z_1$ and $z_2$ be two points on the ellipse $\mathcal{E}$ at which the expression $\vert f(z)-f(z_0)\vert$ attains its maximum and minimum respectively. Then,
\[\lim_{a\to 0}\big\vert \frac{f(z_1)-f(z_0)}{f(z_2)-f(z_0)}\big\vert =1
 \]
\end{itemize}
where $a$ is the great axis of the ellipse.

\end{enumerate}
Lavrentieff says that if we assume that the function $p$ is bounded, then we recover a class of functions which is analogous to the one considered by Gr\"otzsch in \cite{Gr1928}.

After Lavrentieff wrote his Comptes Rendus note, another comptes Rendus note was published by Stoilov, \cite{Stoilov}, showing that the first two conditions in Lavrentieff's  definition of an almost analytic function can be expressed using a notion of \emph{interior transformation} that he had introduced.   Stoilov concludes that for any almost analytic transformation 
\[w=f(z)\]
 in the sense of Lavrentieff, there exists a topological transformation of the domain of the variable, 
\[z=t(z')\]
such that the composed map
\[f(t(z'))\]
is analytic. Thus, from the topological point of view, analytic and almost analytic mappings are the same.

Lavrentieff obtained  the following result stating that from the  two functions $p$ and $\theta$ one can recover the almost analytic function $f$ (Theorem 1 in \cite{Lavrentiev1935a} p. 1011, and Theorem 3, \cite{Lavrentieff1927b} p. 414):
\begin{quote}\small
\emph{Given any two function $p(z)\geq 1$ and $\theta(z)$ defined on the closed unit disc such that $p(z)$ is continuous and $\theta$ is continuous at any point satisfying $p(z)\not=1$, there exists an almost analytic function $w=f(z)$ satisfying $f(0)=0$, $f(1)=1$ and sending homeomorphically the closed unit disc to itself and having the functions $p$ and $z$ as characteristics functions. }
\end{quote}

Lavrentieff obtained the following additional uniqueness theorem (\cite{Lavrentiev1935a} p. 1012), saying that the two functions $p$ and $\theta$ determine $f$ provided these functions coincide on a non-discrete subset:
\begin{quote}
\emph{Two almost analytic functions defined on the same domain $D$ having the same characteristic functions $p$ and $\theta$ and coinciding on a subset which has a limit point in $D$ coincide on the whole set $D$.
}
\end{quote}

The two preceding results are versions of the integrability of almost-complex structures in dimension 2, and they also constitute geometric versions of later existence and uniqueness results for quasiconformal mappings with a given dilatation. A huge literature was dedicated later on to get more precise analytic forms of these results, with a minimum amount of continuity assumptions on the dilatation. These results culminated in the so-called Ahlfors-Bers measurable Riemann mapping theorem and its refinements.

Lavrentieff also obtained  the following result, which has the same flavor as a later result of Teichm\"uller \cite{T200}  in which a growth condition is assumed on the quasiconformal dilatation (Theorem 2 in \cite{Lavrentiev1935a} p. 1011 and Theorem 7 in  \cite{Lavrentiev1935} p. 418):
\begin{quote}\small
\emph{Let $p(z)\geq 1$ and $\theta(z)$  be two functions defined on the closed unit disc such that $p$ is continuous and $\theta$ is continuous at every point satisfying $p(z)\not=1$. For every $r<1$  let $q(r)$ be the maximum of $p(z)$ on the circle $\vert z\vert =r$. If the integral $\int_0^r \frac{dr}{rq(r)}$ is divergent, then one can construct an almost analytic function $w=f(z)$ satisfying $f(0)=0$, $f(1)=1$, whose characteristic functions are $p$ and $\theta$ and which sends homeomorphically the  closed unit disc to itself.}
\end{quote}

There is a condition involving the divergence of the integral $\int_0^r \frac{dr}{rq(r)}$ in another result of  Lavrentieff's \cite{Lavrentiev1935a} which concerns the so-called type problem. We shall quote this result below, but first, we briefly recall the type problem, to which we shall also refer later in the present paper, in our review of the works of Teichm\"uller and of Ahlfors.

 To understand the type problem, one starts by recalling that by the uniformization theorem, every simply connected Riemann surface which is not the sphere is conformally equivalent to either the complex plane or the unit disc. The type problem asks for a way to decide, for a given simply-connected surface defined in some manner (a surface associated to a meromorphic function, defined as a ramified cover, etc.) whether it is conformally equivalent to the complex plane or to the unit disc. In the first case, the surface is said to be of parabolic type, and in the second case of hyperbolic type. In some sense, the type problem is a constructive complement to the uniformization question. There are precise relations between the type problem and Picard's theorem, and we shall mention one of them below.  
  
  Despite its apparent simplicity, the type problem turned out to be very difficult, and it gave rise to a profound theory.  
 The importance of this problem is stressed by Ahlfors in his paper \cite{Ahlfors-JMPA}, in which he formulates the type problem for Riemann surfaces, but also a related ``type problem" for  univalent functions: Such functions are of two types: either infinity is the only point omitted, or it is not.  In the first case, he calls the function \emph{parabolic}, and in the second case he calls it \emph{hyperbolic}. 
He says that instead of talking about parabolic or hyperbolic functions, one could talk about their associated Riemann surfaces, as parabolic and hyperbolic surfaces. He adds that the notion of multiform function (an inverse function of a univalent function is generally multiform) is ``infinitely simpler than that of a Riemann surface." He then writes that ``This problem is, or ought to be, the central problem in the theory of functions. It is evident that its complete solution would give us, at the same time, all the theorems which have a purely qualitative character on meromorphic functions."

Needless to say, to a geometer, the complex plane is ``naturally" Euclidean, and the unit disc is ``naturally" hyperbolic. Making this precise and useful in function theory is another matter, and there are instances where this intuition is misleading, as we shall see e.g. in \S\,\ref{s:T2}, concerning a disproof by Teichm\"uller of a conjecture by Nevanlinna.

  Lavrentieff studied applications of quasiconformal mappings to the type problem in \cite{Lavrentiev1935a} p. 1012 and  \cite{Lavrentiev1935} p. 421. It is possible that this was the first time quasiconformal mappings are used in the type problem. In the 3-dimensional space with coordinates $(x,y,z)$ we consider a  simply connected surface $S$ defined by an equation
$t=t(x,y)$,  
where $x$ and $y$ vary in $\mathbb{R}$ and where $t(x,y)$ is a function of class $C^1$. The question is to find the type of $S$. Lavrentieff proves the following:
 
 \begin{quote}
 \emph{Let $q(r)$ be the maximum of $1+\vert \mathrm{grad}\ t(x,y)\vert$ on the circle $x^2+y^2=r^2$.  If $\int_1^\infty \frac{dr}{rq(r)}$ diverges, then $S$ is of parabolic type.
 }
 \end{quote}
 The proof is based on the existence theorem for almost-analyic functions, Theorem 7 in  \cite{Lavrentiev1935} which we quoted above.
 
 After this result, Lavrentieff gives a construction of surfaces of hyperbolic type:

  \begin{quote}
 \emph{With the above notation, for the surface $S$ to be of hyperbolic type, it suffices that the following condition be satisfied:
 For any domain $D$ contained in $S$ and containing the point $(0,0,t(0,0))$, the area of $D$ is smaller than $l^{2-\epsilon}$ where $l$ is the length of the frontier of $D$  and $\epsilon$ is an arbitrary fixed positive number.}
 \end{quote}

   The paper \cite{Lavrentiev1935} contains further results. Lavrentieff addresses the question of finding a conformal representation of a two-dimensional Riemannian metric defined by a length element of the general form
\begin{equation} \label{eq:R}
ds^2=Edx^2+2Fdxdy+Gdy^2,\ \ EG-F^2>0
\end{equation}
   onto a domain in the Euclidean plane.
   This question is essentially the one that was raised by the ninteenth-century cartographers/differential geometers, which we reviewed in \S \ref{s:qc}. 
   
Lavrentieff assumes that the functions $E,F,G$ are continuous.  He mentions again the 1928 paper \cite{Gr1928}  by Gr\"otzsch, and he declares that they both discovered the same notion of almost analytic functions.

The papers \cite{Lavrentiev1935a} and \cite{Lavrentiev1935} contain other analytic  and geometric applications. In particular, Lavrentieff obtains an analogue of Picard's theorem generalized to almost analytic functions. We recall the statement of Picard's theorem, because quasiconformal mappings play an important role in the later developments of this theorem.

Picard discovered the famous theorem which carries his name in 1879  \cite{Picard1879}; cf. also the paper \cite{Picard1880} in which he gives more details. The result has two parts, the first concerns entire functions (that is, holomorphic functions defined on the whole complex plane), and the second concerns meromorphic functions.  The two results are sometimes called the Picard ``small"  and ``big" theorem  respectively. The result concerning entire functions says that for any non-constant entire function $f$,  the equation $f(z)=a$ has   a solution for any complex number $a$ except possibly for one value. This generalizes in a significant way Liouville's theorem saying that the image of a non-constant entire function is unbounded. The  result concerning meromorphic functions says that such a function, in the neighborhood of an essential singularity (that is, a point where the function has no finite or infinite limit) takes any complex value infinitely often, with at most one exception.

The proofs by Picard of his results use the theory of elliptic functions and results on the Laplace equation.
A large amount of papers were published, starting from 1896 (17 years after Picard obtained his result!), whose aim was to find simpler proofs and generalizations of this theorem. They include Borel \cite{Borel-Picard}, Schottky \cite{Schottky-Picard}, Montel \cite{Montel1912},  \cite{Bloch-Picard}, Carath\'eodory \cite{Caratheodory-Picard}, Landau \cite{Landau1904}, Lindel\"of \cite{Lindelof-Picard}, Milloux \cite{Milloux-Picard}, Valiron \cite{Valiron-Picard}, and there are many others. 
These various geometric proofs and generalizations of Picard's theorem are the basis of the so-called \emph{value distribution theory}, which was born officially in the work of Nevanlinna in the 1920s, at the publication of his papers \cite{N1922}, \cite{N1924} \cite{N4} and his book \emph{Le th\'eor\`eme de Picard-Borel et la th\'eorie des fonctions m\'eromorphes}, \cite{N1}. This theory became one of the most active mathematical theories of the twentieth century. The introduction of quasiconformal mappings in that theory was one of the major steps in its developement.

 Lavrentieff obtained (\cite{Lavrentiev1935a} p. 1012) the following generalization of Picard's theorem:

\begin{quote}\small 
\emph{Let $w=f(z)$ be an almost analytic function defined on the unit disc.  Consider the characteristic functions $p(z)\geq 1$ and $\theta (z)$ and assume that $\theta$ is continuous at every point $z$ satisfying $p(z)\not=1$.  Assume also that $0$ is an essential singularity. Then, in the neighborhood of   $z=0$ the equation $f(z)=a$ has an infinite number of solutions, with the possible exception of one value of $a$.}
\end{quote}
Lavrentieff notes that Gr\"otzsch in \cite{Gr1928} proved the same theorem under the condition that the function $p$ is bounded (this is the strong form of quasiconformality).

There are geometric applications given by Lavrentieff in his paper \cite{Lavrentiev1935}, concerning the Riemann mapping  theorem applied to surfaces equipped with Riemannian metrics given in the form (\ref{eq:R}).

Lavrentieff is also one of the main founders of the theory of quasiconformal mappings in dimensions $\geq 3$, but this is another story.

    The book \cite{Lavrentieff-century} contains very rich information about Lavrentieff, his life and his epoch. It is written by his former students and colleagues. One chapter consists of an edition of Lavrentieff's memories. Lavrentieff writes there (p. 37--38) that in
1929, as a senior engineer at the Theoretical Department of Russian Institute of aerodynamics, he was given the task of determining the velocity field of the fluids, in a problem related to thin wings, and that he wanted in some way to ``justify the mathematics." In a lapse of time of six months, he managed, on the basis of variational principles using conformal mappings, to find a number of estimates for the desired solution. These estimates allowed him to identify a class of functions, among which the solution had to be found. The problem was reduced to the question of 
solving a system of linear equations, and this led to the solution of this problem. It turned out that the theory of conformal mappings could not fully meet the needs of aerodynamics. This was related to the necessity of taking   into account the compressibility of the air and the possibility of exceeding the speed of sound, and this led to the study of a non-linear system of partial differential equations. The theory of conformal mappings needed to be extended to a wider class, and this was the birth of a new theory of quasi-conformal mappings.

In the same book \cite{Lavrentieff-century}, in the chapter by Ibragimova, it is recounted how Lavrentieff,  while he was waiting for his future wife at a tram stop in Riga,  solved a problem with which he had fought unsuccessfully during more than two years; the main point of the solution consisted in the introduction of quasiconformal mappings, in the problem at hand.
 
The reader can find more details on Lavrentieff's work in the papers \cite{Lavrentiev-Sobolev} and \cite{Lavrentiev-Leray}.

\section{Ahlfors}\label{s:A} 

  In this section, after a few words on Ahlfors' life, we highlight a few aspects of his early work on quasiconformal mappings.

   There are several good biographies of Ahlfors and we shall be brief. We refer the interested reader to his short autobiography in Volume I of his \emph{Collected works} edition \cite{Ahlfors-Collected}, to Gehring's biography in the first part of \cite{Gehring-Ahlfors}, to the notes of his work by Gehring, Ossermann, Kra which appeared in the Notices of the AMS \cite{GKKO} and which are reproduced in \cite{Gehring-Ahlfors}, to Lehto's article \cite{Lehto-Int} in the Mathematical Intelligencer and to the lively biography by Lehto which appeard recently in English translation \cite{Lehto-Ahlfors}.

 Ahlfors studied at the University of Helsinki, in the 1920s. This university was already an internationally known center for function theory,  with the presence of Ernst Lindel\"of and Rolf Nevanlinna as professors. Ahlfors obtained his Master degree in the Spring of 1928. In the autumn of the same year, he accompanied Nevanlinna, who was one of his teachers, for a stay in Zurich where the latter was an invited professor in replacement of Hermann Weyl at the \emph{Eidgen\"ossische Technische Hochschule}. Weyl was on leave that year.
  Ahlfors later on described his new environment as follows: ``I found myself suddenly transported from the periphery to the center of Europe" (\cite{Lehto-Int} p. 4). 
  
  In 1929, Ahlfors created a surprise in the mathematical world by proving a 21-year old conjecture by Denjoy concerning the asymptotic values of an entire function. Ahlfors heard about this conjecture at Nevanlinna's lectures in Zurich. The conjecture says that the number of finite asymptotic values of an entire function of order $k$ is at most $2k$. Making sense of this statement needs the definition of an appropriate notion of ``order" and ``asymptotic value" of an entire function, and this is done in the setting of the so-called Nevanlinna theory, or value distribution theory.  Ahlfors' approach to the Denjoy conjecture was completely new, based on quasiconformal mappings. His solution became part of his doctoral thesis, written under Nevanlinna, and submitted in 1930.
In his article ``The joy of function theory" \cite{Ahlfors-Joy}, Ahlfors writes: 
\begin{quote}\small In retrospect, the problem by itself was hardly worthy of the hullabaloo it had caused, but it was not of the kind that attracts talents, especially young talents. It is not unusual that the same mathematical idea will surface, independently, in several places, when the time is ripe. My habits at the time did not include regular checking of the periodicals, and I was not aware that H. Gr\"otzsch had published papers based on ideas similar to mine, which he too could have used to prove the Denjoy conjecture. 
\end{quote}

Ahlfors' work on quasiconformal mappings was first essentially directed towards applications to function theory, which remained his main field of research during more than 10 years (roughly speaking, from 1929 to 1941).  The type problem was one of his central objects of interest during this period.
After that, Ahlfors became thoroughly involved in Teichm\"uller theory and Kleinian groups; this was also the subject of his remarkable collaboration with Bers. 
Regarding this collaboration, Abikoff writes in his paper \cite{Abikoff-Bers}: 
\begin{quote}\small
It was during this period of extensive expository
writing that Bers's research blossomed as well.
In his 1958 address to the International Congress
of Mathematicians, Bers announced a new proof
of the so-called Measurable Riemann Mapping
Theorem. He then essentially listed the theorems
that  followed  directly  from  this  method, including  the  solution  to  Riemann's  problem  of
moduli.  He  was  outlining  the  work  of  several
years, much of which was either joint with or paralleled  by  work  of  Lars  Ahlfors.  Their  professional collaboration was spiritually very close;
they  were  in  constant  contact  although  they
wrote only one paper together, the paper on the
Measurable Riemann Mapping Theorem. Nonetheless,  their  joint  efforts,  usually  independent and often simultaneous, and their personal generosity inspired a cameraderie in the vaguely  defined  group  which  formed  around
them. The group has often been referred to as
the  ``Ahlfors-Bers  family"  or  ``Bers  Mafia";  cf.
the recent article by Kra \cite{Kra-Bers}. It is gratifying to see that the spirit of cooperation they fostered lives on among many members of several generations of mathematicians.
\end{quote}

 Until the end of his life, quasiconformal mappings continued to play an essential role in Ahlfors' work. His 1954 paper \cite{Ahlfors1954},  which concerns a new proof of Teichm\"uller's Theorem, contains several basic results on quasiconformal mappings, including a general definition of such mappings using the distortion of quadrilaterals -- a definition where the differentiability of the map becomes irrelevant.\footnote{Such a definition is also attributed to Pfluger \cite{Pfluger}, independently. It was noticed later that the class of functions obtained by this definition coincides with a class introduced by Morrey in 1938, consisting of weak homeomorphic solutions $f(z)$ of the Beltrami equation $f_{\overline{z}}=\mu(z)f_z$, where $\mu(z)$ is a measurable function satisfying $\sup \vert \mu(z)\vert <1
 $.} The paper also contains a reflection principle, a compactness result, and an analogue of the Hurwitz theorem for quasiconformal mappings.

 Like Lavrentieff did before him, Ahlfors used  quasiconformal mappings as an important tool to prove results in function theory. One of the important by-products of
his 1935 paper paper \emph{Zur Theorie der \"Uberlagerungsfl\"achen} \cite{Ahlfors-Zur} (On the theory of covering surfaces) was to show that several of Nevanlinna's results on value distribution theory had little to do with conformality; in fact they hold for quasiconformal mappings.
This came quite as a surprise, since quasiconformal mappings are very different from the conformal ones. The latter satisfy very strong rigidity properties; for instance, the values of such a mapping in an arbitrary small region determines its values everywhere. Quasiconformal mappings do not satisfy such properties.  Ahflors writes in the introduction to his paper \cite{Ahlfors-Zur} (translation in \cite{Lehto-Ahlfors} p. 28): 
\begin{quote}\small 
This work had its origin in my endeavor to get by geometric means the most significant results of meromorphic function theory. In these attempts it became evident that only an easily limited portion of R. Nevanlinna's Main Theorems, and thereby nearly all classical results, were dependent on the analyticity of the mapping. In contrast, their entire structure is determined by the metric and topological properties of the Riemann surface, which is the image of the complex plane. The image surface is then thought of being spread out above the Riemann sphere, i.e. to be the covering surface of a closed surface.
\end{quote}

Ahlfors also introduced the techniques of quasiconformal geometry to prove results on conformal function theory. The need for a space of more flexible functions to work with was one reason for which quasiconformal mappings naturally appeared in the theory of functions of a  complex variable.

Another important point of view which Ahlfors brought in function theory is topology. The topological techniques of branched coverings, cutting and pasting pieces of the complex plane,  the Euler characteristic, the Riemann-Hurwitz formula, and other topological tools, became with Ahlfors part of the theory of functions of a complex variable. In some sense, Ahlfors' new approach constituted a return to the sources, that is, to the geometrical and topological methods introduced by Riemann, for whom meromorphic functions were  nothing else than branched coverings of the sphere characterized by some finite data describing their singularities. Carath\'eodory, in presenting Ahlfors' work at the Fields medal ceremony in 1936, declared that Ahlfors opened up a new chapter in analysis which could be called ``metric topology" (quoted in \cite{Lehto-Int} p. 4). The theorems of Picard  and Nevanlinna on value distribution became, under Ahlfors, theorems on ``islands": instead of counting the number of times that certain values are omitted, in Ahlfors' theory, one counts how many ``regions" are omitted.  
Ahlfors also brought the Gauss-Bonnet theorem and other results of differential geometry into the realm of Nevanlinna's theory.
 
Let us finally mention that Ahlfors' first works on the type problem (see \cite{Ahlfors1931} and \cite{Ahlfors-CRAS}), unlike the works of Laverentieff on the same problem, involves the so-called length-area method which is a common tool used by Gr\"otzsch and himself, and later on by Teichm\"uller, in relation with quasiconformal mappings. Regarding this method, Ahlfors writes, commenting on his first two published papers on asymptotic values of entire functions of finite order in his \emph{Collected Works} \cite{Ahlfors-Collected} Vol. 1 (p. 1): 
\begin{quote}
The salient feature of the proof is the use of what is now called the length-area method. The early history of this method is obscure, but I knew it from and was inspired by its application in the well-known textbook of Hurwitz-Courant to the boundary correspondence in conformal mapping. None of us [Ahlfors is talking about Nevanlinna and himself] was aware that only months earlier H. Gr\"otzsch had published two important papers on extremal problems in conformal mapping in which the same method is used in a more sophisticated manner. My only priority, if I can claim one, is to have used the method on a problem that is not originally stated in terms of conformal mapping. The method that Gr\"otzsch and I used is a precursor of the method of extremal length [...]  In my thesis \cite{Ahlfors:Unter} the lemma on conformal mapping has become the main theorem in the form of a strong and explicit inequality or distortion theorem for the conformal mapping from a general strip domain to a parallel strip, together with a weaker inequality in the opposite direction. [...] A more precise form of the first inequality was later given by O. Teichm\"uller.
\end{quote}
Ahlfors refers here to Teichm\"uller's paper \cite{T200}, on which we comment in what follows.

\section{On Teichm\"uller's writings} \label{s:T}
       
During his short lifetime, Oswald Teichm\"uller (1913--1943) published a series of fundamental papers on geometric function theory and on moduli of Riemann surfaces. At the same time, he developed to a high degree of sophistication the theory of quasiconformal mappings and their applications. He also wrote several papers on algebra and number theory. His writings were read by very few mathematicians. His article \emph{Extremale quasikonforme Abbildungen und quadratische Differentiale} (Extremal  quasiconformal mappings and quadratic differentials)  \cite{T20}, published in 1939, is known among Teichm\"uller theorists, but it is seldom quoted in the mathematical literature. In fact, to this day, mathscinet mentions only 35 citations of that paper: 6 between 1949 and 1984, and the rest (29 citations) between 2000 and 2016. These figures are ridiculously low, in view to the major impact that the paper had and the extensive literature to which it gave rise. Furthermore, when Teichm\"uller's papers were quoted, they were generally accompanied by a little bit of suspicion, and the impression which resulted from the comments on them is that Teichm\"uller did not provide proofs for his results. Bers, in his 1968 ICM address \cite{Bers1960A}, referring to the paper \cite{T20} we just quoted, writes the following:   
\begin{quote} \small
Much of [my] work consists in clarifying and verifying assertions of Teichm\"uller whose bold ideas, though sometimes stated awkwardly and without complete proofs, influenced all recent investigations, as well as the work of Kodaira and Spencer on the higher dimensional case.
\end{quote}  Ahlfors, in his 1954 paper \cite{Ahlfors1954}, writes: 
\begin{quote} \small
In a systematic way the problem of extremal quasiconformal mapping was taken up by Teichm\"uller in a brilliant and unconventional paper \cite{T20}. He formulates the general problem and, although unable to give a binding proof, is led by heuristic arguments to a highly elegant conjectured solution. The paper contains numerous fundamental applications which clearly show the importance of the problem.
\end{quote}

It is a matter of fact that Teichm\"uller cared much more about transmitting ideas than about writing detailed proofs. In the paper \cite{T20}, he introduced the space that became later known as \emph{Teichm\"uller space}, he equipped it with the metric that bears the name \emph{Teichm\"uller metric} and he thoroughly investigated its  properties, including its infinitesimal (Finsler) structure, and he proved the existence and uniqueness of geodesics between any two points. He developed the theory of quasiconformal mappings between arbitrary Riemann surfaces and the partial differential equations they satisfy,  giving a characterization of the tangent and cotangent spaces at each point of Teichm\"uller space as a space of equivalence classes of partial differential equations and quadratic differentials respectively, he also studied  Teichm\"uller discs (which he called complex geodesics), he investigated global convexity properties of the space. There are many other ideas and results in that paper, including the development of quasiconformal invariants of Riemann surfaces. We shall comment on these ideas in \S \ref{s:T2} below.

There are other remarkable papers of Teichm\"uller, and all of them are almost never quoted. Among the most important ones, we mention the paper \emph{Untersuchungen \"uber konforme und quasikonforme Abbildungen} (Investigations on conformal and quasiconformal mappings) \cite{T200} in which Teichm\"uller brings  new techniques and results on quasiconformal geometry with applications to Nevanlinna's theory, the paper \emph{Ver\"anderliche Riemannsche Fl\"achen} (Variable Riemann surfaces) \cite{T32}, in which he lays down the bases of the complex structure of Teichm\"uller space, and the paper  \emph{\"Uber Extremalprobleme der konformen Geometrie} (On extremal problems in conformal geometry)
 \cite{T23} which contains a broad research program in which Teichm\"uller displays analogies between Riemann surface theory, algebra and Galois theory. The exact content of all these papers and the methods developed there are very poorly known, even among the specialists.  Our aim in the next section is to give an overview of some of Teichm\"uller's important ideas on quasiconformal mappings and their applications in geometric function theory, Riemann surfaces and moduli. Many of these papers are still worth reading in detail today; they contain ideas that lead directly to interesting research projects. 
 Let us now say a few words on Teichm\"uller's writings.

In the first 40 years after their publication, Teichm\"uller's works attracted almost exclusively the interest of analysts, represented in a masterly manner by Ahlfors and Bers. 
Geometers became interested in Teichm\"uller's work essentially when Thurston entered the scene, in the 1970s. 
Let us quote Thurston, from his  foreword to Hubbard's book on Teichm\"uller spaces \cite{Hubbard}:
\begin{quote}\small
Teichm\"uller theory is an amazing subject, richly connected to geometry, topology, dynamics, analysis and algebra. I did not know this at the beginning of my career: as a topologist, I started out thinking of Teichm\"uller theory as an obscure branch of analysis irrelevant to my interests. My first encounter with Teichm\"uller theory was from the side. I was interested in some questions about isotopy classes of homeomorphisms of surfaces, and after struggling for quite a while, I finally proved a classification theorem for surface homeomorphisms. [...] Bers gave a new proof of my classification theorem by a method that was much simpler than my own, modulo principles of Teichm\"uller theory that had been developed decades earlier. 
From this encounter I came to appreciate the beauty of Teichm\"uller theory, and of the close connections between 1-dimensional complex analysis and two and three-dimensional geometry and topology. A great deal of mathematics has been developed since that time and there are many active connections between geometry, topology, dynamics and Teichm\"uller theory.
\end{quote}

Teichm\"uller's papers have a conversational and informal style, but the details of the proofs are difficult to follow. At several places, difficult statements are given with only sketches of proofs. At some points, the proofs, when they exist, involve new ideas from various directions including geometry, topology, partial differential equations, etc. This makes them arduous for a reader whose knowledge is specialized. Teichm\"uller sometimes states a problem using several equivalent formulations, explains the difficulties and  the methods of attack, returns to the same problem several pages later, sometimes several times, each time from a new point of view and with a new explanation, and eventually comes up with a new idea and declares that the problem is solved, after he had presented it as a conjecture. In the introduction of the paper  \emph{Bestimmung der extremalen quasikonformen Abbildungen bei geschlossenen orientierten Riemannschen Fl\"achen}
(Determination of extremal quasiconformal mappings of closed oriented Riemann surfaces) \cite{T29}  (1943), he writes, referring to his previous paper \emph{Extremale quasikonforme Abbildungen und quadratische Differentiale} \cite{T20}:   
\begin{quote}\small
In 1939, it was a risk to publish a lengthy article entirely built on conjectures. I had studied the topic thoroughly, was convinced of the truth of my conjectures and I did not want to keep back from the public the beautiful connections and perspectives that I had arrived at. Moreover, I wanted to encourage attempts for proofs. I acknowledge the reproaches, that have been made to me from various sides, even today, as justifiable but only in the sense that an {\it unscrupulous\/} imitation of my procedure would certainly lead to a barbarization of our mathematical literature. But I never had any doubts about the correctness of my article, and I am glad now to be able to actually prove its main part.

At that time, I did not have an exact theory of {\it modules\/}, the conformal invariants of closed Riemann surfaces and similar ``principal domains." In the meantime I have developed such a theory
aiming at the intended application to quasiconformal mappings. I will have to briefly report on it elsewhere. The present proof does not depend on this new theory, and  works with the notion of {\it uniformization\/} instead. However, I think one will have to combine both to bring the full content of my article \cite{T20} in mathematically exact form.
\end{quote}
At the end of the same paper, he writes:
\begin{quote}\small
[...] I cannot elaborate on these and similar questions  and only express my opinion that all those aspects that have been treated separately in the discussion so far will appear  as a great unified theory of variable Riemann surfaces in the near future.
\end{quote}
Teichm\"uller died soon after, without being able to develop the unified theory he was aiming for. 

     Ahlfors and Gehring write in the preface to Teichm\"uller's \emph{Collected works}:
\begin{quote}\small
Teichm\"uller's style was unorthodox, to say the least. He himself was well aware of the difference between a proof and an intuitive reasoning, but his manner of presentation makes it difficult to follow the frequent shifts from one mode to another.
\end{quote}
Abikoff, in his expostory article on Teichm\"uller \cite{abikoff}, recalls that the aim of \textit{Deutsche Mathematik}, the journal  in which most of Teichm\"uller's papers were published, was to concentrate on general ideas and not on technical details. This was supposed to be the tradition to which the name ``German mathematics" refers, a tradition that was shared by Riemann and Klein. Abikoff adds that he learned from a conversation with Herbert Busemann that  ``Teichm\"uller manifested those traits early in his career but when pressed could offer a formal proof." In fact, at some places, Teichm\"uller, with his style, exceeded the norms of \emph{Deutsche Mathematik}. For instance, the editors of that journal included the following footnote in the introduction of the paper \emph{ \"Uber Extremalprobleme der konformen Geometrie}\cite{T23} (1941):  
  \begin{quote}\small 
  The following work is obviously unfinished, it has the character of a fragment.
     Unreasonably high demands are made on the reader's cooperation and imagination.
     For the assertions, which are not even stated precisely with rigor,
 neither proofs nor even any clues are given.
 One thing not of fundamental importance,  the ``residual elements'' occupy
a broad space for something almost unintelligible, while far too scarce hints are provided for
 fundamentally important individual examples.
  But the author explains that he cannot do better in the foreseeable future. --
 If we still publish the   author's remarks, despite all lack that distinguishes the work against the other papers in this journal,
it is to bring up for discussion the thoughts contained therein relating to the theme of estimates for schlicht functions.
\end{quote}
         
   Concerning the purely mathematical content and as a general rule, Teichm\"uller avoided the use of sophisticated tools, and his proofs are based on first principles. This is probably another reason for which his papers are difficult to read. Quoting theorems and building upon them often makes things easier than developing a theory from scratch. Teichm\"uller writes, in the paper  \cite{T29} (1943):  
        \begin{quote}\small
        It is rather because of my unaptness and lack of knowledge on continuous topology, that I can only prove the purely topological lemma taking such long detours. But in my opinion it would be absurd to introduce the length of the required proof as a measure for the value of a statement. Of course none of the necessary details of a proof  may be left out if a statement seems worth a proof at all. But as our example shows clearly, there is a natural division in main points and minor points that we may not arbitrarily change. 
\end{quote}

           Let us now say a few words about Teichm\"uller's life.
           
Teichm\"uller entered the university of G\"ottingen in 1931, at a time where this university was at its full glory. Hilbert had retired from his professorship the year before, but was still lecturing  on philosophy of science. Hermann Weyl, who came back from Zurich,  officially as Hilbert's successor, was among Teichm\"uller's teachers. He was lecturing on differential geometry, algebraic topology and  the philosophy of mathematics.  Teichm\"uller's other teachers included Richard Courant (who was the head of the institute), Otto Neugebauer, Gustav Herglotz,  Emmy Noether and Edmund Landau.  Besides the professors, there were brilliant assistants and graduate students at the mathematics department, including Saunders Mac Lane, Herbert Busemann, Ernst Witt, Hans Wittich, Werner Fenchel and Fritz John. Max Born and James Franck were among the professors at the physics department. 

The Nazis came to power in Germany in February 1933 and the deluge arrived very rapidly. On April 7, 1933, a law was voted, excluding the Jews and the ``politically unreliable" from civil service. In particular, Jewish faculty members at all state universities had to be immediately dismissed, with a few exceptions  (those who served the German army during World War I and those who had more than 25 years of academic service). At the same time, a state law restricted the number of Jewish students at German universities. The University of G\"ottingen immediately declined. It continued though to receive distinguished visitors. In particular, Rolf Nevanlinna lectured there in
1936 and 1937. Teichm\"uller was very 
much influenced by these lectures. 

Teichm\"uller obtained his doctorate two years later, in 1935, under Helmut Hasse. The subject was functional analysis. Right after his thesis, he published several papers on algebra and geometric algebra. In October 1936, he started working on his Habilitation under Bieberbach, and he obtained it in 1938. The subject was function theory. The habilitation is published as \cite{T200}.
Most of Teichm\"uller's papers on function theory and Riemann surfaces were published during the Second World War (between 1939 and 1943). In the summer of 1939,  and just after he was offered a position at the University of Berlin, he was drafted into the army. The war soon broke out and he was sent to the Norwegian front. He spent there a few weeks, after which he returned to Berlin, but he never officially took a job at the university. Until the year 1943, he occupied a position at the high command of the army. 

It seems however that he was not satisfied with this work at the army. In 1941, he applied for a position at Munich's  Ludwig-Maximilians University, to take over the position which would eventually to  vacant after Carath\'eodory's retirement. The paper \cite{Litten1994} by Freddy Litten describes the complex events related to Carath\'eodory's  succession. The book \cite{Georgiadou} by Maria Georgiadou also contains several pages devoted to this question. This author writes (p. 355) that ``in the six years preceding the appointment of his successor, Carath\'eodory was active in proposing or blocking candidacies." At a certain stage, Helmuth Kneser and Teichm\"uller were on the top of the list. We read in \cite{Georgiadou} p. 360 (information obtained from \cite{Litten1994}, note 166 p. 163) that Wilhelm F\"uhrer, who was the mathematics referee at the education ministry regarding that position, described Teichm\"uller as an ``old Nazi, free from any kind of authoritarian behavior, and Kneser as a party comrade and an active SA member." We also read in \cite{Georgiadou} p. 361 (information from \cite{Litten1994}, note 166 p. 163ff.)  that ``Carath\'eodory complained about the imperfect reasoning in a treatise of Teichm\"uller's." On November 13, 1941, Carath\'eodory wrote to Wilhelm S\"uss (a former Nazi and member of the SA, who was the president of the DMV -- the German Mathematical Society -- and who by the way is the founder of the Oberwolfach Mathematics Research Institute) that ``Teichm\"uller was not sufficiently mature for one of the main posts at a great German university" (\cite{Georgiadou} p. 361, from Nachla\ss \ Wilhelm S\"uss C89/48, Freiburg Universit\"atsarchiv). Despite Carath\'eodory's reservations, at a certain stage, Teichm\"uller was ranked first. But at a later senate's session, ``it was revealed that the head of the lecturers had altered documents in Teichm\"uller's favor, so Teichm\"uller was crossed off the list" (\cite{Georgiadou} p. 361, information from \cite{Litten1994}, note 166 p. 165ff.). Thus, Teichm\"uller did not get the position at Munich. In the fall of 1943, he volunteered as a soldier at the Eastern Front, and he presumably died there, the same year. In a set of notes, written by Teichm\"uller's mother and translated by Abikoff in \cite{abikoff},  one can read: ``Even as a soldier he found the time to write papers while his comrades rested." 
In the introduction of his paper \cite{T23}, Teichm\"ulle writes: 
\begin{quote}\small
Because I only have a limited vacation time at my disposal, I cannot give reasons for many things, but only assert them.
This is unfortunate because I do not know anyway
  the exact and generally valid proofs of the principles to be drawn up.
After all, the expert familiar with \emph{Extremale quasikonforme Abbildungen und quadratische Differentiale}
will be able to add much of what is missing. Also I  have not yet been able to perform  many individual studies.
\end{quote}

                            In the introduction of \cite{T20}, he writes: 
                            
                            \begin{quote}\small
                            As I said before, the main results will not be proven in the strictest sense.
I only hope to establish them in a way that any serious doubts are practically excluded, and to encourage finding the proofs.
Thus when possible, the ideas that led to discovering the solution are conveyed in direct chain.
As far as I can see, attempts of exact proofs find their right place only where this train of thought ends, according to the paradox ``proving is reversing the train of thought."
So it is expected that a systematic theory progressing toward exact proofs would cause \emph{for the time being} bigger difficulties in understanding than the present heuristic introduction.
\end{quote}

Many people were reluctant to read Teichm\"uller's papers, because they had heard of his  extreme political ideals.  
Concerning his active involvement in the Nazi movement during his years in G\"ottingen, there are several versions and explanations. Abikoff collected information from Busemann, Fenchel and others who were present at G\"ottingen at that period. In his paper \cite{abikoff}, he reports that Fenchel recalls: ``That Teichm\"uller was a member of the Nazi party, we learned when he distributed Nazi propaganda in the Mathematics Institute.  Otto Neugebauer who assisted Courant in the administration of the Institute threw him out."  In a recent book \cite{Lehto-Ahlfors}, Lehto writes (p. 74): ``Hans Wittich had been in G\"ottingen at the same time as Teichm\"uller. Both had listened to Rolf Nevanlinna's lectures concerning Riemann surfaces and had had mathematical collaboration. Wittich who had many contacts in Z\"urich, never joined the National Socialist party, and soon after the Allied forces deemed him politically fit to serve as a professor. In his view, Teichm\"uller could not have had many chances to be politically active during his time at G\"ottingen since he spent all days sitting in the library with his mathematical papers." Lehto adds: ``But many dissenting opinions were voiced, particularly among Jewish mathematicians. They acknowledged Teichm\"uller's mathematical merits but regarded him as a human being of the basest sort possible. In various ways they aimed to show that Teichm\"uller's mathematical achievements were not as groundbreaking as was given to believe since numerous mathematicians before him had worked on similar questions."

In conclusion, there are several reasons for which results of Teichm\"uller remained unknown for a long time. Some of them still remain so. As a consequence, some of these results were rediscovered by others. They are contained in papers written by experts in the field, without any mention of Teichm\"uller's name and probably without any knowledge of the fact that they had already been discovered by him.\footnote{It is a fact of experience that in the mathematical  literature, references and original sources are often uncontrolled by the authors. As soon as a  mathematician gives a reference, it is likely that several other authors will repeat the reference without checking the original paper. Thus, an error in an attribution usually propagates, especially if the real author of the result is no more alive, and even if his name is Euler or Poincar\'e.
 The author of the present article fell a few times into this trap. For instance, in \cite{Handbook-Intro}, he wrote, repeating what is usually claimed in the literature, that ``the bases of the complex analytic theory of Teichm\"uller space were developed by Ahlfors and Bers" without any mention of Teichm\"uller, whose paper \cite{T32} is entirely dedicated to this subject and lays the bases of this theory, many years before Ahlfors and Bers worked on the subject.}  It is important to make straight the history of science and to attribute the origin of an idea to the person who discovered it first.

\section{Teichm\"uller's work on quasiconformal mappings} \label{s:T2}
 
It is not our intention  in the rest of this paper to write a detailed summary of Teichm\"uller's results nor their impact, and we shall leave aside many theorems he proved. Instead, we concentrate on some of the ideas he brought in the theory of quasiconformal mappings and their applications. The reader who wishes to learn more about Teichm\"uller's specific results is referred to his original papers. Some of them exist now in English translation and appeared in various volumes of the Handbook of Teichm\"uller theory   \cite{T20,T23, T24, T29, T31,T32, T200}.  Others will appear in the same series. We also refer to the commentaries \cite{T20C,T29C,T32C,T23C,T31C,T24C, T200C} which appeared in the same volumes.

 We have divided our survey of Teichm\"uller's contribution to quasiconformal mappings into the following subsections:
 \begin{enumerate}
 \item The notion of quasiconformal invariant and the behavior of conformal invariants under quasiconformal mappings.
 \item The existence and uniqueness theorem for extremal quasiconformal mappings.
 
 \item The theory of infinitesimal quasiconformal mappings.
  \item Quasiconformal mappings in function theory.
  
 \item The non-reduced quasiconformal theory.
 \item A distance function on Riemann surfaces using quasiconformal mappings.
 
 \item Another problem on quasiconformal mappings

 \item A general theory of extremal mappings motivated by the extremal problem for quasiconformal mappings.
  \end{enumerate}
 
We now elaborate on these topics.

\subsection{The notion of quasiconformal invariants and the behavior of conformal invariants under quasiconformal mappings.}

In \cite{T20}, Teichm\"uller considers a general surface of finite type, oriented or not, with or without boundary, and with or without distinguished points in the interior or on the boundary. He denotes by $\mathfrak R^{\sigma}$ the Riemann moduli space of that surface; $\sigma$ is the dimension of the space. He was aware of the fact that $\mathfrak R^{\sigma}$ is an orbifold (and not a manifold). He  introduces a covering, Teichm\"uller space, which he shows is a manifold, and which he denotes by $R^{\sigma}$, and he gives the following formula for the dimension: 
$$
\sigma=-6+6g+3\gamma+3n+2h+k+\rho
$$
where  
\begin{itemize}
\item[] $g$ is the number of handles;
\item[] $\gamma$ is the number of crosscaps;
\item[] $n$ is the number of boundary curves;
\item[] $h$ is the number of distinguished points in the interior;
\item[] $k$ is the number of distinguished points on the boundary;
\item[]  $\rho$ is the dimension of the group of conformal self-mappings of the Riemann surface.
\end{itemize}

There are several relations with quasiconformal geometry, which we recall now.

 In the introduction of the paper \cite{T20},  Teichm\"uller asks in all generality and in a precise manner the question of finding mappings between two arbitrary surfaces that are closest to conformal, in relation with the general question of the behavior of \emph{conformal invariants}. He writes:
\begin{quote} \small
In the present study, \emph{the behavior of conformal invariants under quasiconformal mappings} shall be examined.
This will lead to the problem of \emph{finding the mappings that deviate as little as possible from conformality under certain additional conditions}.
\end{quote} 
 
Teichm\"uller was able to develop these ideas after he introduced   Teichm\"uller space, equipping it with a topology.  We note in passing that Riemann's cryptic statement about the existence of $3g-3$ moduli was a parameter count, and can hardly be considered as a dimension count, since there was no topology and no manifold structure, let alone complex structure within sight. We refer the interested reader to the surveys \cite{2015a} and \cite{ji&papadop} on the birth of the topology and the complex structure on Teichm\"uller space.
 
 We briefly recall some facts on the topology that Teichm\"uller introduced on Teichm\"uller space, and the resulting quotient topology on Riemann's moduli space.  This is no closely related to the subject of our survey, since these topologies are based on the notion of quasiconformal mapping.
 
  The topology is induced by the so-called Teichm\"uller metric, where the distance between two  Riemann surfaces is the logarithm of the least dilatation of all quasiconformal mappings between them.  This is defined in \S 18 of the paper \cite{T20}.   (Note that the metric, at this point of the paper, is defined on moduli space, and not yet on Teichm\"uller space.)  We shall recall the precise statement below. In \S 19, Teichm\"uller formulates the following questions in which $P$ and $Q$ denote points in moduli space:
 
\begin{quote}
\emph{What is the distance $[PQ]$ between $P$ and $Q$?
In particular, can we determine all extremal quasiconformal mappings whose dilatation quotients are everywhere $\leq e^{[PQ]}$?}

It is expected that there always exists an extremal quasiconformal mapping and that it is unique up to conformal self-mappings of the given principal regions.
\end{quote}
This is the main ``conjecture" which Teichm\"uller eventually proved. The result became known under the name Teichm\"uller existence and uniqueness theorem. He obtained the proof in the paper \cite{T20} and in the later paper \cite{T29}.
    
 Teichm\"uller
introduces a \emph{conformal invariant} as a function on Riemann's moduli space $\mathfrak R^{\sigma}$ (\S 13).  He says that locally, there are precisely $\sigma$ independent conformal invariants. In \S 16, he talks about the behavior of conformal invariants under quasiconformal mappings. From the definition, a conformal invariant undergoes little changes under quasiconformal mappings if the supremum of the dilatation quotient is sufficiently close to 1. Teichm\"uller asks the following: 

 \begin{quote}
{\it Let $J$ be a conformal invariant, considered as a function on $\mathfrak R^{\sigma}$. Consider a given Riemann surface, and let $C>1$ be a given real number.
What are the values that $J$ takes at those Riemann surfaces onto which the given Riemann surface can be quasiconformally mapped so that the dilatation quotient remains everywhere $\leq C$?}
 
 \end{quote}

This conceptually important question is a formulation of one of the main problems he was interested in. At the end of \S 17, he   says that the notion of quasiconformal dilatation can be used to \emph{define} the topology of Riemann's moduli space. This is explained more precisely in \S 18, where Teichm\"uller introduces his distance. He writes:

\begin{quote}
 
{\it Let us consider two Riemann surfaces of the same type represented in $\mathfrak R^{\sigma}$ by the points $P$ and $Q$.
Given any quasiconformal mapping between these Riemann surfaces, let $C$ denote the supremum of its dilatation quotient.
We define the distance $[PQ]$ between the two points, or between the two Riemann surfaces, as the logarithm of the infimum of all these $C$.}

 \end{quote}

Teichm\"uller formulates a more general question (end of \S 16) where the quasiconformal bound $C$ of a quasiconformal mapping is not a real number, but a \emph{function} on the surface. He says that this general question is hard. Such a theory of quasiconformal mappings with dilatation  not necessarily bounded uniformly was developed later on by various authors, e.g. Lehto in \cite{L1970} and \cite{L1971}, David in \cite{David1988} and Brakalova and Jenkins in \cite{BJ1}; see also the survey by Otal in Volume IV of the Handbook of Teichm\"uller theory  \cite{Otal}. A notion of quasiconformal mappings whose dilatation is not bounded by a constant but by a function with a specific behavior at infinity is also used in Teichm\"uller's paper \cite{T200}. Here, he proves that if a one-to-one map from the complex plane to itself has a dilatation quotient $D(z)$ which satisfies $D(z)\leq C(\vert z\vert)$ where $C(r)$ is a real function defined for $r\geq 0$ satisfying $C(r)\to 1$ as $r\to \infty$ such that 
\[\int^\infty (C(r)-1)\frac{dr}{r}<\infty,\] then as $z\to \infty$ we have 
\[\vert w\vert \sim  \mathrm{const} \cdot \vert z\vert.\]
We noted that Lavrentieff, before Teichm\"uller,  in his paper \cite{Lavrentiev1935} obtained a related result. Lavrentieff's definition of quasiconformal mappings did not involve a uniform bound  of the dilatation. 

 In his study of the behavior of conformal invariants under quasiconformal mappings, Teichm\"uller a result on hyperbolic lengths of closed curves. In \S 35 of \cite{T29}, he  establishes an inequality which compares the length of a simple closed geodesic on a hyperbolic surface to the length of the corresponding geodesic when the hyperbolic surface  is transformed into another one by a quasiconformal mapping. The inequality involves the dilatation of the quasiconformal mapping. This inequality was rediscovered by Sorvali in \cite{sorvali} (1973) and then by Wolpert in \cite{wolpert} (1979). The result is known today as \emph{Wolpert's inequality}. Teichm\"uller obtains this inequality using the point of view of Fuchsian groups; he analyses the effect of a quasiconformal mapping conjugating two Fuchsian groups on the dilatations of the hyperbolic transformations. The same method was used by Sorvali and Wolpert.

The point of view on the behavior of conformal invariants under quasiconformal mappings allows a variety of new definitions of quasiconformal mappings. Indeed, one starts with a conformal invariant of the surface (a Green function, extremal length, hyperbolic length, etc.) One then declares that a homeomorphism of the surface is quasiconformal if it distorts the given conformal invariant by a bounded amount. In the 1967 paper \cite{Gehring1967} by Gehring, based on this idea, 12 definitions of quasiconformal mappings are given.

Let us note, to end this section on conformal invariants, that in \S 17 of \cite{T20}, Teichm\"uller considers a problem related to surfaces of infinite type. He says that it is possible to start with the conformal invariants of a family of Riemann surfaces of finite type and study the behavior of limits of these invariants on families of Riemann surfaces that converge to a Riemann surface of infinite type.

\medskip

\subsection{The existence and uniqueness theorem for extremal quasiconformal mappings.}

In his paper \emph{Bestimmung der extremalen quasikonformen Abbildungen bei geschlossenen orientierten Riemannschen Fl\"achen} (Determination of extremal quasiconformal mappings of closed oriented Riemann surfaces) \cite{T29} (1943), Teichm\"uller completed the proof of his major result concerning extremal quasiconformal mappings. The theoremsays that in any homotopy class of homeomorphism between two closed orientable Riemann surfaces of genus $g\geq 2$ there exists a unique extremal quasiconformal mapping, and that this mapping, in appropriate local coordinates, is affine, that is, it has the form:
\[
z \mapsto K\cdot\mathrm{Re}\left(z \right) +\mathrm{Im}\left( z \right).\]

The local coordinates on the surface are given by the horizontal and the vertical measured foliations of a holomorphic quadratic differential associated to the homeomorphism. The resulting extremal mapping is usually called the \emph{Teichm\"uller mapping}. As we already noted, the result was stated  in the paper \cite{T20} as a ``conjecture" and in the setting of an arbitrary surface of finite type (with or without boundary, and with or without distinguished points in the interior and on the boundary).  The uniqueness part of the statement was proved in that paper. The existence was  obtained for closed surfaces in \cite{T29}. In the same paper, Teichm\"uller says that he will publish later on a proof of his theorem for an arbitrary surface of finite type.

One of the main tools in Teichm\"uller's existence and uniqueness theorem is the notion of  holomorphic quadratic differential.  The \emph{complex dilatation} of an extremal quasiconformal mapping is expressed in terms of such a differential.  in \S 46 of his paper, Teichm\"uller begins a study of the connection between quasiconformal mappings and quadratic differentials. He reports the following:
\begin{quote}\small

During a night of 1938, I came up with the following {\it conjecture}:

{\it Let $d\zeta^{2}$ be an everywhere finite quadratic differential on $\frak{F}$, different from $0$.
Assign to every point of $\frak{F}$ the direction where $d\zeta^{2}$ is positive.
Extremal quasiconformal mappings are described through the direction fields obtained this way and through arbitrary constant dilatation quotients $K\geq1$.}

\end{quote}
This dream that he had was the starting point of an extremely rich theory which transformed the field of quasiconformal theory.

 The appearance of holomorphic quadratic differentials to describe the dilatation of a quasiconformal mapping  is in itself surprising because a quasiconformal mapping is a non-holomorphic object.  Ahlfors writes in his paper \cite{Ahlfors-Joy}, p. 445:
 \begin{quote}\small
 Over the years this important notion [quasiconformal mappings] has changed the nature of function theory quite radically, and so on many different levels. The ultimate miracle was performed by O. Teichm\"uller, a completely unbelievable phenomenon for better and for worse. He managed to show that a simple extremal problem, which deals with quasiconformal mappings of a Riemann surface, but does not involve any analyticity, has its solution in terms of a special class of analytic functions, the quadratic differentials. Over the years since Teichm\"uller's demise his legacy has mushroomed to a new branch of mathematics, known as Teichm\"uller theory, whose connection with function theory is now almost unrecognizable.
 \end{quote}
 On the other hand, a quadratic differential on a surface is characterized by a pair of orthogonal measured foliations. It is good to recall here that a pair of orthogonal foliations on a Riemann surface, in relation with extremal problems for mappings between surfaces, already appear in the work of Tissot; see \S\ref{s:Tissot}.
 
In the theory of quadratic differentials that Teichm\"uller developed, the consideration of distinguished points is a significant issue. Technically speaking, this involves the study of the behavior of quadratic differentials at these points, which is an important aspect. Such a theory was further developed by Jenkins in several papers. We shall mention some of these works below. Teichm\"uller's idea led to an outline of a far-reaching theory of extremal maps, with relations with algebra and Galois theory which Teichm\"uller outlined in his paper \cite{T23} and on which we shall also report below.

Teichm\"uller's existence and uniqueness theorem generalizes Gr\"otzsch's existence and uniqueness result of a quasiconformal mapping with least dilatation between quadrilaterals. Teichm\"uller, in \cite{T24}, gave a proof of an analogous theorem for the case of a pentagon, using the so-called method of continuity. In the last section of the paper  \cite{T29}, he writes: ``During my research, I have always kept in mind the aim to give a continuity proof for the existence of quasiconformal mappings similar to the one in the case of the pentagon." In his 1964 survey paper on quasiconformal mappings and their applications \cite{Ahlfors1964}, Ahlfors, writes that the result on pentagons is ``already a sophisticated result.'' Teichm\"uller writes (\S 1 of his paper \cite{T24}) that the case of the pentagon is ``the simplest case of the higher cases,'' and that this simple case already ``shows how far one has to go beyond and extend the methods of Ahlfors and Gr\"otzsch.''

Teichm\"uller's existence and uniqueness theorem for arbitrary closed surfaces of genus $g\geq 2$ has far-reaching applications. For instance, it shows that in the Teichm\"uller metric, any two points are connected by a unique geodesic and it gives a proof that the Teichm\"uller space of a closed surface of genus $g$ is homeomorphic to $\mathbb{R}^{6g-6}$. Note however that Teichm\"uller, in his paper \cite{T29}, uses the dimension count, which he had already obtained in \cite{T20} using other methods, to prove his existence theorem.  

Ahlfors and Bers tried, during many years, to give simpler and more accessible proofs of Teichm\"uller's theorem. 
In his paper \cite{Bers1960}, written 18 years after Teichm\"uller's paper,  Bers writes:
\begin{quote}\small
Our arrangement of the arguments preserves the logical structures of Teichm\"uller's proof; the details are carried out differently. More precisely, we work with the most general definition of quasiconformality, we rely on the theory of partial differential equations in some crucial parts of the argument, and we make use of a simple set of moduli for marked Riemann surfaces.
\end{quote}

 In another paper, \cite{B2}, written 16 years after the paper \cite{Bers1960} which we just mentioned, and which is among the very last papers that Bers wrote, Bers declares again that all the arguments given by Teichm\"uller in his paper \cite{T29} are correct.   In the abstract of that paper, he says that he gives a ``simplified version of Teichm\"uller's proof (independent of the theory of Beltrami equations with measurable coefficients) of a proposition underlying his continuity argument for the existence part of his theorem on extremal quasiconformal mappings." He then writes, in the introduction: 
 \begin{quote}\small
 In proving a basic continuity assertion (Lemma 1 in \S\,14C of \cite{Bers1960}) I made use of a property of quasiconformal mappings (stated for the first time in \cite{B0} and also in \S\,4F of \cite{Bers1960}) which belongs to the theory of quasiconformal mappings with bounded measurable Beltrami coefficients (and seems not to have been known to Teichm\"uller). Some readers concluded that the use of that theory was indispensable for the proof of Teichm\"uller's theorem. This is not so, and Teichm\"uller's own argument is correct.
\end{quote}

         In Volume 2, p. 1 of the edition of his \emph{Collected works} \cite{Ahlfors-Collected}, commenting on his  paper  \cite{Ahlfors1954},  Ahlfors writes: 
         \begin{quote}\small
         More than a decade had passed since Teichm\"uller wrote his remarkable paper \cite{T20} on extremal quasiconformal mappings and quadratic differentials. It has become increasingly evident that Teichm\"uller's ideas would profoundly influence analysis and especially the theory of functions of a complex variable, although nobody at that time could foresee the extent to which this would be true. The foundations of the theory were not commensurate with the loftiness of Teichm\"uller's vision, and I thought it was time to reexamine the basic concepts. My paper has serious shortcomings, but it has nevertheless been very influential and has led to a resurgence of interest in quasiconformal mappings and Teichm\"uller theory. 
Based on this definition the first four chapters are a careful and rather detailed discussion of the basic properties of quasiconformal mappings to the extent that they were known at that time. In particular a complete proof of the uniqueness part of Teichm\"uller's theorem was included. Like all other known proofs of the uniqueness it was modeled on Teichm\"uller's own proof, which used uniformization and the length-area method. Where Teichm\"uller was sketchy I tried to be more precise.
In the original paper Teichm\"uller did not prove the existence part of his theorem, but in a following paper \cite{T29} he gave a proof based on a continuity method. I found his proof rather hard to read and although I did not doubt its validity I thought that a direct variational proof would be preferable. My attempted proof on these lines had a flaw, and even my subsequent correction does not convince me today. In any case my attempt was too complicated and did not deserve to succeed.
\end{quote} 

All this gives an idea of  the fate of one of the results of Teichm\"uller's paper \cite{T20}, and of Ahlfors and Bers's involvement in the task of making Teichm\"uller's results clearer, several decades after they were written.

   \subsection{The theory of infinitesimal quasiconformal mappings.}

In \S 54 of his paper  \cite{T20},  Teichm\"uller introduces the notion of {\it infinitesimally quasiconformal mapping}, and he develops this theory thoroughly in the later sections of the paper. Instead of being a mapping, an infinitesimally quasiconformal mapping is a vector field on the surface.  Teichm\"uller develops the relation between these objects and the so-called Beltrami differentials. He works with conformal structures underlying Riemannian metrics, written locally as quadratic forms in the tradition of Gauss, and he studies the effect of quasiconformal metrics on the coefficients of these quadratic forms. In fact, he develops the partial differential equation approach to quasiconformal mappings. In particular, he determines the expression of the metric after an infinitesimal perturbation by an infinitesimally quasiconformal mapping. The deformation is expressed in terms of a function $B$ which depends on local coordinates, and Teichm\"uller notes that the expression $B \frac{dz^2}{\vert dz \vert^2}$ does not depend on them. Later in the paper, $B$ will be considered as a tangent vector to Teichm\"uller space, that is, an element of a vector space he denotes by $L^\sigma$ and which will play an extremely important role. Teichm\"uller defines a norm $\|\cdot\|$  on the space $L^\sigma$. The extremal problem is stated in a new form, using the notion of extremal infinitesimal quasiconformal mapping.
Teichm\"uller announces that the solution of this problem, in the case of closed surfaces, will be given in terms of a holomorphic quadratic differential. This is done in the subsequent sections of the paper.

  In \S 68 of the same paper, Teichm\"uller introduces the notion of \emph{locally extremal} differentials. He highlights the link between extremal quasiconformal mappings and the orthogonal pair of foliations on a surface that are given by the quadratic differential. He then defines an equivalence relation on the set of $B\frac{dz^2}{\vert dz \vert^2}$, and a \textit{pairing} between these forms and quadratic differentials. He associates to each class of $B\frac{dz^2}{\vert dz \vert^2}$, $3g-3$ complex numbers $\left( k_{\mu} \right)_{1\leq \mu \leq \tau}$. He showed that conversely, for  any complex numbers $k_1 , k_2 , \cdots , k_\tau$, there exists one and only one class of invariant forms $B\frac{dz^2}{\vert dz\vert^2}$ having the desired properties. He deduces that the space $L^\sigma$ is a complex vector space of dimension $3g-3$. He asks whether the norm $\|\cdot\|$  comes from a Hermitian structure. This question was answered by Weil, who introduced the so-called Weil-Petersson structure on Teichm\"uller space.
 We refer the reader to the paper \cite{T20C} which is a commentary on these questions.
 
 Teichm\"uller's theory of infinitesimal quasiconformal mappings is at the basis of the theory of deformation of higher-dimensional complex structures that was developed later on By Kodaira and Spencer.

\subsection{Quasiconformal mappings in function theory.}

The applications of the theory of quasiconformal mappings  to function theory are highlighted by Teichm\"uller since the beginning of his paper \emph{Extremale quasikonforme Abbildungen und quadratische Differentiale} \cite{T20} which we mentioned several times. In the introduction, he writes: 
 
 \begin{quote} \small
 
 As I have done it before, I shall here also examine quasiconformal mappings not exclusively for their own sake, but chiefly because of their connections with notions and questions that interest function theorists. It is true that quasiconformal mappings have been for only very few years systematically applied to purely function theoretical questions, and so this method has managed to gain until now only a limited number of friends.
It may nevertheless be mentioned that the contribution I have been able to make some time ago in this journal to the Bieberbach coefficients problem relies in fact on ideas that shall be developed here.\footnote{[Teichm\"uller's footnote] O. Teichm\"uller, Ungleichungen zwischen den Koeffizienten schlichter Funktionen, these proceedings, 1938. (This is Item \cite{T14} in the bibliography of the present paper.)}

 \end{quote}
 
In the last sections of the same paper (\S 167--170), Teichm\"uller reviews several applications of quasiconformal mappings to function theory.   He also says, in relation with these applications, that one should not only consider quasiconformal mappings  with bounded dilatation quotient,  but also those satisfying some appropriate upper estimates of the dilatation quotient (that is, this quotient is bounded pointwise,  by some function defined on the surface). We already commented on this point of view.

 The paper \emph{Eine Anwendung quasikonformer Abbildungen auf das Typenproblem} (1937) \cite{T9} is one of the first papers that Teichm\"uller wrote on function theory, and it concerns the type problem.\footnote{Drasin notes in \cite{Drasin1986} that Teichm\"uller was the first to realize the use of quasi-conformal mappings in the type problem. This remark is not correct, because as we noted, Lavrentieff did this before Teichm\"uller.} He proves there that the type of a simply connected Riemann surface is invariant by quasiconformal mapping.  This is based on the fact that there is no quasiconformal mapping between the unit disc and the complex plane. It seems that this fact is proved for the first time in the paper \cite{T9}.

Teichm\"uller also introduces in the paper  \cite{T9} the theory of quasiconformal mappings in the study of line complexes. A line complex is a combinatorial graph on a Riemann surface which is given in the form of an infinite branched covering of the sphere, with a finite number of branching points. Such an object was introduced by Nevanlinna in \cite{N-Ueber}  and Elfving in \cite{Elfving1934}, and it has the same flavor as an object introduced by Speiser in \cite{Speiser1930} and which became known as a \emph{Speiser tree}. The line complex captures the combinatorics of the branched covering, but it does not determine uniquely the surface. The reader might remember that the description of a Riemann surface as a branched covering of the sphere, provided with the combinatorial information at the branch points, is in the pure tradition of Riemann, in particular, of the work done in his doctoral thesis \cite{Riemann-Grundlagen}.  We briefly recall the definition of a line complex. 

We start with a simply connected surface $\frak{M}$ given as an infinite branched covering of the Riemann sphere. We let $a_1,\ldots,a_q$ be the set of branching values (this is a subset of the Riemann sphere). We assume that this set of branching values is finite.

We let $L$ be a simple closed curve on the sphere that passes once through each of the branching values points $a_1,\ldots,a_q$. The sphere is then decomposed into two subsets $\frak{I}$ and $\frak{A}$ which are homeomorphic to closed discs. The covering surface $\frak{M}$ is obtained by gluing copies of simply-connected pieces called ``half-planes," each one covering either the disc $\frak{I}$ or the disc $\frak{A}$. Thus, a half-plane is a copy of  either $\frak{I}$ or $\frak{A}$. The Riemann surface $\frak{M}$ is assumed to be an infinite branched covering of the sphere, but the number of half-planes that are glued among each other at a given branching point of $\frak{M}$ may be finite or infinite. The line complex is a graph embedded in  $\frak{M}$ which encodes the overall gluing of half-planes that constitutes this surface $\frak{M}$. To define this graph, we choose a point in the interior of each of the two polygons $\frak{I}$ and $\frak{A}$; let these points be $P_1$ and $P_2$ respectively.   Let $s_1,\ldots,s_q$ be the $q$ arcs that are cut on the closed curve $L$ by the  points $a_1,\ldots,a_q$. We join $P_1$ and $P_2$ be $q$ simple arcs, each of them crossing the interior of one of the arcs  $s_i$ ($i=1,\ldots, q$) transversely and no other point of $L$.  We obtain a graph on the Riemann sphere, which has two vertices and $q$ edges connecting them. The inverse image of this graph by the covering map from $\frak{M}$ to the sphere is the line complex associated to the covering.  

The line complex does not completely determine the surface $\frak{M}$, and the problem which Teichm\"uller addresses in the paper \cite{T9} is whether one can determine, from the line complex associated to a the surface  $\frak{M}$, the type of  $\frak{M}$. He says that this problem is still far from being solved. His main contribution in this paper is the fact that any two such surfaces $\frak{M}$ and $\frak{M}'$ that have the same line complex have the the same type. The proof of this result is done in two steps, first showing that two surfaces $\frak{M}$ and $\frak{M}'$ with the same line complex are mapped quasiconformally to each other and then that the complex plane and the unit dic are not quasiconformally equivalent.

The paper \emph{Untersuchungen \"uber konforme und quasikonforme Abbildungen}  (Investigations of conformal and quasiconformal mappings) \cite{T200} (1938) is also on function theory. Teichm\"uller gives a criterion for a class of surfaces to be of hyperbolic type, and a negative response to a conjecture by Nevanlinna concerning line complexes. A review of this result is contained in  Chapter II of Nevanlinna's book \emph{Analytic functions} (p. 312), where the latter says that Teichm\"uller answered in the negative a question he raised, asking ``whether a transcendental simply connected surface described as a certain infinite ramified covering is parabolic or hyperbolic according to whether its angle geometry is Euclidean or Lobachevskian." We give a rough idea of the ideas that are involved in this work because they are interesting for low-dimensional geometers and topologists. 

           Nevanlinna's conjecture is based on the fact that the type of the simply connected surface $\frak{M}$ (we use the above notation), given a branched covering of the sphere,  is related to the degree of ramification of its line complex. If the ramification is small, the surface will be of parabolic type, and if it is large, it will be of hyperbolic type. The motivation behind this conjecture is that a large degree of ramification will put a large angle structure at the vertices of the graph $\Gamma$, which corresponds to negative combinatorial curvature. Nevanlinna writes (p. 309 of \cite{N-analytic}) that ``it is natural to imagine the existence of a critical degree of ramification that separates the more weakly ramified parabolic surfaces from the more strongly ramified hyperbolic surfaces." Indeed, Nevanlinna introduced a notion of combinatorial curvature of a line complex, which he called the \emph{excess}. His conjecture (\cite{N-analytic} p. 312) says that the surface $\frak{M}$ is parabolic or hyperbolic according as the mean excess of its line complex is zero or negative.   Teichm\"ulller in his paper \cite{T200}, disproved the conjecture by exhibiting a simply connected surface $S$ ramified over the sphere whose corresponding line complex is hyperbolic and whose mean excess is zero.

Teichm\"uller's paper \emph{Ungleichungen zwischen den Koeffizienten schlichter} \cite{T14} (1938) is related to the so-called coefficient problem, or Bieberbach conjecture. In this paper, Teichm\"uller initiated the use of quadratic differentials in the study of the coefficient problem. In the same paper, he obtained a sequence of inequalities which were supposed to prove the Bieberbach conjecture. Abikoff, in \cite{abikoff} writes the following:  ``Teichm\"uller \cite{T14} gave a heuristic argument for the existence of the quadratic differential in the course of proving a somewhat bizarre collection of inequalities satisfied by the coefficients $a_n$.  He claims that a proper use of these inequalities leads to the solution of Bieberbach's conjecture; however neither he nor any other user of this technique ever came close to succeeding."

\subsection{The non-reduced quasiconformal theory.}
The \emph{non-reduced quasiconformal Teichm\"uller theory} is the Teichm\"uller theory that includes the  quasiconformal techniques that are involved in the setting of surfaces with boundary in which every point on the boundary is considered as a distinguished points. In particular, the homotopies which define the equivalence relation between Riemann surfaces (marked by quasiconformal mappings) fix pointwise the boundary. An example of a non-reduced Teichm\"uller space is the Teichm\"uller space of the disc, which is in essence the so-called \emph{universal Teichm\"uller space}.  Such a theory was already conceived by Teichm\"uller. In fact, in his paper \cite{T20}, Teichm\"uller mentions an even more general theory, namely, he considers arbitrary bordered  Riemann surfaces with distinguished arcs on the boundary where every point of such an arc is considered as a distinguished point. The later paper \cite{T31}, which we discuss below, also deals with a problem that belongs to this setting.   

 In the non-reduced theory, there is a  \emph{non-reduced Teichm\"uller space}, and a \emph{non-reduced moduli space}. There are nontrivial problems that arise in the non-reduced theory which have no analogue in the reduced theory. For instance, there was a long-standing conjecture about the general form of extremal mappings in the non-reduced theory.  Teichm\"uller suggested in \cite{T20} (p 185) that the extremal maps of the disk, with some given boundary condition, are Teichm\"uler mappings.  V. Bo\v{z}in, N. Laki\'c, V. Markovi\'c and M. Mateljevi\'c showed in \cite{yougos} that there exist extremal mappings (and even, uniquely extremal ones) which are not Teichm\"uller mappings and which are associated to quadratic differentials which may be of finite or infinite norm. These authors obtained conditions under which the extremal mapping is unique \cite{yougos}  \cite{yougos1}  \cite{MM}.

\subsection{A distance function on Riemann surfaces using quasiconformal mappings.}

   In \S 27 of the paper \cite{T20}, Teichm\"uller  defines a distance function on a sphere with three distinguished points which  is based on the notion of quasiconformal mappings. Precisely, the distance between two points is the logarithm of the least quasiconformal constant of a mapping of the surface that carries one point to another. He shows that this distance coincides with the hyperbolic distance. In \S 160, he gives an analogous definition in the case of an arbitrary closed surface, and he states that up to a condition, the surface, equipped with this distance, is a Finsler manifold.
   
 Kra rediscovered this distance, and  studied it in his paper \cite{kra} in which he showed that the Bers fibres are not Teichm\"uller discs. 
Kra proved that an extremal mapping that realizes the distance between two points is unique if these points are close enough for the hyperbolic metric. He also showed that the distance is equivalent to the hyperbolic metric but not proportional to it unless the surface is the thrice-punctured sphere. The distance is called the Kra distance in \cite{Shen}.

\subsection{Another problem on quasiconformal mappings}
  
In the paper 
\emph{Ein Verschiebungssatz der quasikonformen Abbildung (A displacement theorem of  quasiconformal mapping)}   \cite{T31}, published in 1944, Teichm\"uller solves a problem on extremal quasiconformal mappings that has many ramifications and which is in the spirit of ideas contained in \S27 of \cite{T20} which we just mentioned.
The question is to find an extremal quasiconformal self-map of the disk which satisfies the following two properties:
 \begin{itemize}
 \item  the image of 0 is $-x$, for some prescribed $0<x<1$;
  \item the restriction of the map to the unit circle is the identity.
  \end{itemize}
  Teichm\"uller describes an extremal mapping for this problem. He shows that this mapping  is unique and that it is a Teichm\"uller map. It is associated to a quadratic differential which has finite norm and a singularity which is a simple pole.
  
  We recall that a conformal mapping of the disc which fixes pointwise the boundary is necessarily the identity. In the present setting, one relaxes the conformality property, and asks for the best quasiconformal mapping which is the identity on the boundary and which satisfies a certain property for some points in the interior.
  This is an extremal problem in the setting the non-reduced theory. The existence of the extremal quasiconformal mapping is reduced to the problem of obtaining an affine map between the surfaces bounded by ellipses.  The solution of the problem is again described in terms of a quadratic differential. The dilatation of the extremal mapping is shown to be constant on the surface. In \cite{Reich1987}, Reich gave a generalized formulation of this problem, replacing, from the beginning,  the unit disc by the interior of an ellipse. He solved a more general problem concerning mappings of conics. In the case where the conic is a parabola, this leads to a case where the extremal mapping is non-unique.
  Ahlfors, in his papers \cite{Ahlfors-Bounded} and  \cite{Ahlfors-Open}, considered similar extremal problems.

Teichm\"uller's  paper is commented on by Vincent Alberge in \cite{T31C}.
 \subsection{A general theory of extremal mappings motivated by the extremal problem for quasiconformal mappings.}

 Teichm\"uller's paper \emph{\"Uber Extremalprobleme der konformen Geometrie (On extremal problems in conformal geometry)} \cite{T23}, published in 1941, contains ideas that generalize some of the ideas expressed in Teichm\"uller's paper \cite{T20} on the use of quadratic differentials to solve the extremal problem of quasiconformal mappings. His aim is to apply his techniques to general geometric extremal problems. 
 Ahlfors made strong statements on the importance of general extremal problems. He writes in \cite{Ahlfors-dev} p. 500:
\begin{quote}\small
I have frequently mentioned extremal problems in conformal mapping, and I believe their importance cannot be overestimated. It is evident that extremal mappings must be the cornerstone in any theory which tries to classify conformal mappings according to invariant properties.
\end{quote}

In \cite{Ahlfors-SIAM} p. 3, Ahlors writes:
\begin{quote}\small
In complex function theory, as in many other branches of analysis, one of the most powerful classical methods has been to formulate, solve, and analyze extremal problems. This remains the most valuable tool even today, and constitutes a direct link with the classical tradition.
\end{quote}

 The generality of the extremal mappings that Teichm\"uller alludes to in his paper \cite{T23} is much beyond that of the extremal mappings to which Ahlfors refers.
  
Teichm\"uller claims in  \cite{T23} that some of his ideas expressed in his previous papers  on moduli of Riemann surfaces and on algebra are related to each other, that they may lead to  general concepts, and that they are applicable to the coefficient problem for univalent functions. These idea involves the introduction of a new structure at the distinguished points of a surface, viz., an ``order'' for the series expansions of functions in the local coordinate charts at these distinguished points (the Riemann surface is characterized, as Riemann showed, by the field of meromorphic function it carries). The general principle is that if at the distinguished point  the extremal problem requires a function that has a fixed value for its first  $n$ derivatives, then the quadratic differential should have a pole of order $n+1$ at that point. This kind of idea  was used by Jenkins in his approach to the Bieberbach conjecture, cf. \cite{J1, J2, J3, J5}; we discuss this fact again below.

 The paper \cite{T23} is difficult, and it was probably never read thoroughly by any other mathematician than its author. Teichm\"uller starts by commenting on the fact that function theory is closely related to topology and algebra. For instance, one is led, in dealing with function-theoretic questions, to prove new generalizations of the Riemann-Roch theorem. The emphasis, he says, should not be on the Riemann surface but rather on the marked points. He makes an analogy with a situation in algebra, that he had considered in \cite{T4}, namely, one is given three objects, $\mathfrak{A}, \mathfrak{A}_1,\mathfrak{A}_2$; in the geometric case,  $\mathfrak{A}$ is a Riemann surface with distinguished points, $\mathfrak{A}_1$ is the support of the Riemann surface, that is, the Riemann surface without points distinguished, and $\mathfrak{A}_2$ the set of distinguished points. In the algebraic context,  $\mathfrak{A}$ is a normal (Galois) extension of a field,  $\mathfrak{A}_1$ is a cyclic  field extension and $\mathfrak{A}_2$ is the set of generators of the Galois group.  At a distinguished point on a Riemann surface, one chooses a local coordinate $z$ and assumes that the other local coordinates $\tilde{z}$ are of the form
  \[\tilde{z}=z+a_{m+1}z^{m+1}+a_{m+2}z^{m+2}+\ldots
\]
  where $m$ is an integer. This puts a restriction on the first $m$ derivatives of a function belonging to the field of meromorphic functions associated to the Riemann surface. Te definition is done so that this restriction does not depend on the choice of the   local coordinates. Such a distinguished point is said to be of \emph{order} $m$.  
  
Teichm\"uller notes that the results apply to abstract function fields instead of Riemann surfaces, and that estimates on the coefficients of a univalent function may be obtained through a method involving extremal mappings associated with quadratic differentials with some prescribed poles.  He makes relations with several classical problems, including the question of finding the Koebe domain of a family of holomorphic functions defined on the disk, that is, the largest domain contained in the image of every function in the family.

It is possible that his main interest in developing this theory arose from the Bieberbach conjecture.\footnote{We recall that Teichm\"uller wrote his habilitation under the guidance of Bieberbach.}  We recall that this conjecture concerns the coefficients of holomorphic injective functions defined on the unit disc $D=\{z\in \mathbb{C} \vert \ \vert z\vert <1\}$ by a Taylor series expansion: 
 \[f(z)=\sum_{n=0}^\infty a_nz^n
  \]
 normalized by $a_0=0$ and $a_1=1$. The conjecture was formulated by Bieberbach in 1916 and proved fully by Louis de Branges in 1984. It says that the coefficients of such a series satisfy the inequalities
 \[\vert a_n\vert \leq n\]
 for any $n\geq 2$.
 Bieberbach proved in his paper \cite{BB} the case $n=2$. He also showed that equality is attained for the functions of the form $K_\theta (z)=z/(1-e^{i\theta}z)^2$ with $\theta\in\mathbb{R}$. Such a function if the composition of the so-called ``Koebe function''  
 \[k(z)=z/(1-z)^2=z+2z^2+3z^3+\ldots
 \]
with a rotation. In a footnote (\cite{BB} p. 946), Bieberbach suggested that the value $n=\vert a_n(K_\theta)\vert$ might be an upper bound for all the functions satisfying the given assumptions.  
 
 Bieberbach's result is related to the so-called ``area theorem,'' which is referred to by Teichm\"uller in the paper \cite{T23}, a  theorem which gives the so-called ``Koebe quarter theorem" saying that the image of any univalent function $f$ from the unit disc of $\mathbb{C}$ onto a subset of $\mathbb{C}$ contains the disc of center $f\left( 0\right)$ and radius $\vert f^{\prime}\left( 0 \right)\vert /4$. The important fact for the subject of our present article is that for this geometric form of the Bieberbach conjecture, Teichm\"uller introduced the techniques of extremal quasiconformal mappings and quadratic differentials. This is the content of the last part of this paper.
  
   Jenkins, in his series of papers \cite{J1, J2, J3, J5} and others, developed an approach to the coefficient theorem unsing quadratic differentials and based on the works of Teichm\"uller and Gr\"otzsch.  In a 1962 ICM talk \cite{J3}, Jenkins writes: 
   \begin{quote}\small
   Teichm\"uller enunciated the intuitive principle that the solution of a certain type of extremal problem for univalent functions is determined by a quadratic differential for which the following prescriptions hold. If the competing mappings are to have a certain fixed point the quadratic differential will have a simple pole there. If in addition fixed values are required for the first $n$ derivatives of competing functions the quadratic differential will have a pole of order  $n+1$ at the point. He proved a coefficient result which represents a quite special case of the principle but did not obtain any general result of this type. The General Coefficient Theorem was presented originally as an explicit embodiment of Teichm\"uller's principle, that is, the competing functions were subjected to the normalizations implied by the above statement.
   \end{quote}

In \cite{J3}, Jenkins states a theorem which gives a more precise form of a result stated by Teichm\"uller. This result concerns a Riemann surface of finite type equipped with a quadratic differential, with a  decomposition of the surface into subdomains defined by the trajectory structure of this differential. There is a mapping from each of these subdomains onto non-overlapping subdomains of the surface, and these mappings are subject to conditions on preservation of  certain coefficients at the poles. He obtains an inequality that involves the coefficients of the quadratic differential at poles of order greater than one and those of the mapping, with a condition for the inequality to be an equality. This condition states that the function must be an isometry for the metric induced on the surface by the quadratic differential. It is followed by a detailed analysis of the equality case.  Jenkins also refers to numerous specific applications of such a result which Teichm\"uller has conjectured, see \cite{J1}, p. 278-279.

 After the  work done by Teichm\"uller, the theory of quasiconformal mappings exploded and became one of the major tools in function theory, in terms of the volume of the work done, of the difficulties and the deepness of the results, and of its applications in the other fields of mathematics.
 In the survey that we made of the work of Teichm\"uller on quasiconformal mappings, we tried to convey the fact that several of his ideas may still be further explored, and that reading his work will be beneficial.

 \noindent {\bf Acknowledgements.} I would like to thank Vincent Alberge and Olivier Guichard who read a preliminary version of this paper, Mikhail M. Lavrentiev, Jr. for a correspondence about his grandfather, and Nikolai Abrosimov for his help in reading documents in Russian about Lavrentiev.

\printindex
\end{document}